\documentclass[12pt]{amsart}
\usepackage{amssymb,amsmath,amsfonts,latexsym}
\usepackage{enumerate}
\usepackage{hyperref, xcolor}

\newcommand{\newsection}[1]{\setcounter{equation}{0} \section{#1}}
\newcommand{\bea}{\begin{eqnarray}}
\newcommand{\eea}{\end{eqnarray}}

\newcommand{\clb}{\mathcal{B}}

\newcommand{\cld}{\mathcal{D}}

\newcommand{\clh}{\mathcal{H}}

\newcommand{\clm}{\mathcal{M}}
\newcommand{\cln}{\mathcal{N}}

\newcommand{\cls}{\mathcal{S}}

\newcommand{\D}{\mathbb{D}}
\newcommand{\N}{\mathbb{N}}
\newcommand{\C}{\mathbb{C}}

\def\textmatrix#1&#2\\#3&#4\\{\bigl({#1 \atop #3}\ {#2 \atop #4}\bigr)}
\def\dispmatrix#1&#2\\#3&#4\\{\left({#1 \atop #3}\ {#2 \atop #4}\right)}
\newcommand{\be}{\begin{equation}}
\newcommand{\ee}{\end{equation}}
\newcommand{\ben}{\begin{eqnarray*}}
\newcommand{\een}{\end{eqnarray*}}

\newcommand{\NI}{\noindent}

\newcommand{\bi}{\begin{itemize}}
\newcommand{\ei}{\end{itemize}}
\newcommand{\diag}{\mbox{diag}}

\newcommand\la{\langle}
\newcommand\ra{\rangle}

\newtheorem{Theorem}{\sc Theorem}[section]
\newtheorem{Lemma}[Theorem]{\sc Lemma}
\newtheorem{Proposition}[Theorem]{\sc Proposition}
\newtheorem{Corollary}[Theorem]{\sc Corollary}
\newtheorem{Definition}[Theorem]{\sc Definition}
\newtheorem{Question}{\sc Question}
\newtheorem{ass}[Theorem]{\sc Assumption}

\theoremstyle{definition}
\newtheorem{Example}[Theorem]{\sc Example}
\newtheorem{Remark}[Theorem]{\sc Remark}

\newtheorem{Note}[Theorem]{\sc Note}
\newcommand{\bt}{\begin{Theorem}}
\def\beginlem{\begin{Lemma}}
\def\beginprop{\begin{Proposition}}
\def\begincor{\begin{Corollary}}
\def\begindef{\begin{Definition}}
\def\beginexamp{\begin{Example}}
\def\beginrem{\begin{Remark}}
\def\beginq{\begin{Question}}
\def\beginass{\begin{ass}}
\def\beginnote{\begin{Note}}
\newcommand{\et}{\end{Theorem}}
\def\endlem{\end{Lemma}}
\def\endprop{\end{Proposition}}
\def\endcor{\end{Corollary}}
\def\enddef{\end{Definition}}
\def\endexamp{\end{Example}}
\def\endrem{\end{Remark}}
\def\endq{\end{Question}}
\def\endass{\end{ass}}
\def\endnote{\end{Note}}

\numberwithin{equation}{section}

\begin{document}

\title[Invariant subspaces of some compressions of the Hardy shift]{Invariant subspaces of compressions of the Hardy shift on some parametric spaces}

\author[Das]{Susmita Das}
\address{Indian Institute of Science, Department of Mathematics, Bangalore, 560012,
India}
\email{susmitadas@iisc.ac.in, susmita.das.puremath@gmail.com}
\date{\today}


\begin{abstract}
We study the class of operators $S_{\alpha,\beta}$ obtained by compressing the Hardy shift on the parametric spaces $H^2_{\alpha, \beta}$ corresponding to the pair $\{\alpha,\beta\}$ satisfying $|\alpha|^2+|\beta|^2=1$.  We show, for nonzero $\alpha,\beta$, each $S_{\alpha,\beta}$ is indeed a shift $M_z$ on some analytic reproducing kernel Hilbert space and present a complete classification of their invariant subspaces. While all such invariant subspaces $\clm$ are cyclic, we show, unlike other classical shifts, they may not be generated by their corresponding wandering subspaces $(\clm\ominus S_{\alpha,\beta}\clm)$. We provide a necessary and sufficient condition along this line and show, for a certain class of $\alpha, \beta$, there exist $S_{\alpha,\beta}$-invariant subspaces $\clm$ such that $\clm\neq [\clm\ominus S_{\alpha,\beta}\clm]_{S_{\alpha,\beta}}$.
\end{abstract}

\maketitle

\newsection{Introduction}\label{sec: intro}
A common theme in the study of invariant subspaces of shift operators on classical function spaces is the relation between wandering subspaces for the operator in question and the associated invariant subspaces. The starting point of this investigation is, of course, the celebrated theorem of A. Beurling (\cite{Beurling}). To begin with, let us recall the definition of a wandering subspace. Throughout this paper, we assume all the Hilbert spaces $\clh$ are infinite-dimensional, separable and over the complex field $\mathbb{C}$. Following Halmos (\cite{Halmos}), we call a subspace $\cln$ of a Hilbert space $\clh$ is wandering for an operator $T$ if $\cln\perp T^n\cln$ for all $n\in\mathbb{N}$. Therefore, if a subspace $\clm\subseteq\clh$ is $T$-invariant, then $(\clm\ominus T\clm)$ is always wandering for $T$. Corresponding to any subset $E\subseteq\clh$, $[E]_T$ will denote the smallest closed $T$-invariant subspace generated by $\{T^nf:  f\in E, n=0,1,2,\ldots\}$, and we say, a $T$-invariant subspace $\clm$ satisfies wandering subspace property (w.s.p. for short) if $\clm=[\clm\ominus T\clm]_T$. Again, if a wandering subspace $\cln$ satisfies $[\cln]_T=\clh$, we call $\cln$ a generating wandering subspace for $T$. Note that, if $\cln$ is a generating wandering subspace and $\cln=\mathbb{C}f$ for some $f(\neq 0)\in\clh$, $\clh$ is cyclic. With this terminology, let us now state Beurling's Theorem:
\begin{Theorem}
Let $\clm\neq\{0\}$ be a closed subspace of $H^2(\D)$ invariant under the unilateral shift $S$. Then $\dim(\clm\ominus S\clm)=1$ and there exists an inner function $\theta$ such that $$\clm=[\theta]_{S}=[\clm\ominus S\clm]_S.$$
\end{Theorem}
Hence, \textsf{ the map $\clm\longrightarrow \clm\ominus S\clm$ sets a one-one correspondence between the invariant and wandering subspaces of $S$ on $H^2(\D)$. Also, all nonzero wandering subspaces are one-dimensional and generated by an inner function.} Beurling's Theorem plays a crucial role in Operator and Complex Function Theory. It is only natural to look for analogues of the above result in other function spaces. We mention two pertinent examples. First is the Dirichlet space. In (\cite{Richter}), it is shown that, the shift invariant subspaces on Dirichlet space $\cld$, satisfy the similar properties like Hardy-shift i.e., \textsf{$M_z$-invariant subspaces on $\cld$ are generated by their corresponding wandering subspaces and the nonzero wandering subspaces are of dimension one.} Second is the Bergmann space. Unlike the case of Hardy and Dirichlet spaces, the shift $M_z$ on Bergmann space $L^2_a(\D)$ (\cite{Aronszajn}) exhibits different behaviour. The wandering subspaces of $M_z$ on $L^2_a(\D)$ can be of any dimension. In fact, it is shown in (\cite{ABFP}): \textsf{For every $n\in\mathbb{N}\bigcup\{\infty\}$, there exists an invariant subspace $\clm\subseteq L^2_a(\D)$ such that $\dim(\clm\ominus z\clm)=n$.} Also, if an $M_z$-invariant subspace $\clm$ in $L^2_a(\D)$ is cyclic, then $\dim(\clm\ominus z\clm)=1$. Hence, not every Bergmann shift invariant subspaces are cyclic. However, despite this contrast with Hardy and Dirichlet shifts, any $M_z$-invariant subspace on $L^2_a(\D)$ is still generated by the corresponding wandering subspace. This result is proved in (\cite{Aleman-Richter-Sundberg}) and formally states as:
\begin{Theorem}
Let $\clm\subseteq L^2_a(\D)$ be $M_z$-invariant. Then $\clm=[\clm\ominus z\clm]_{M_z}.$
\end{Theorem}
Therefore, in all these three cases of classical shift operators on Hardy, Dirichlet and Bergmann spaces, we note that: \textsf{ A shift-invariant subspace $\clm$ can be recovered from its wandering subspace $\big(\clm\ominus z\clm$\big).} At this point, one may ask:

\NI\textit{What if we perturb the classical shifts a bit? If the resulting operator ( obtained by perturbing a classical shift) is a shift on some function Hilbert space, what can be said about their invariant and wandering subspaces? How (if at all) are these two subspaces related? }
Here shift refers to an operator unitarily equivalent to the multiplication operator $M_z$ on some analytic reproducing kernel Hilbert space. For the detailed theory on reproducing kernel Hilbert spaces, we refer (\cite{Aronszajn}). Throughout, $S$ will denote the usual unilateral shift on $H^2(\D)$. Our aim in this paper is to seek analogues of the above results for the compressions of the unilateral shift $S$ on the closed subspaces $H^2_{\alpha,\beta}$ of the Hardy space. We observe that, these compressions are indeed shifts and can also be represented as perturbations of the Hardy shift. Given a pair of complex number $\alpha, \beta$ with $|\alpha|^2+|\beta|^2=1$, the space $H^2_{\alpha, \beta}$ is a closed subspace of the Hardy space $H^2(\D)$ defined as $$H^2_{\alpha,\beta}=\overline{\text{span}}\{\alpha+\beta z, z^n\quad n\geq 2\}.$$ These spaces are studied in (\cite{Paulsen Raghupati}), in the context of constrained Nevanlinna-Pick interpolation problem. If $\alpha\neq 0$, it is easy to verify that $H^2_{\alpha, \beta}$ is not invariant under the multiplication by $H^{\infty}$-functions. However if $\alpha=0$, then $H^2_{\alpha, \beta}=H^2_{0,1}=z H^2(\D)$, which is invariant under the multiplication by $H^{\infty}$. It is well-known (\cite{Paulsen Raghupati}): for $\alpha\neq 0$, the set of multipliers of $H^2_{\alpha, \beta}$ is precisely $H^{\infty}_1$, where $$H^{\infty}_1= \{f\in H^{\infty}: f'(0)=0\},$$ and also, the spaces $H^2_{\alpha,\beta}$ serve as models for the spaces that are modules over $H^{\infty}_1$. \textsf{Throughout this paper, we will assume both $\alpha, \beta$ are nonzero.} The spaces $H^2_{\alpha,\beta}$ corresponding to such $\{\alpha,\beta\}$ are not invariant under the multiplication $M_z$ of the Hardy space $H^2(\D)$. We are interested in studying the class of operators obtained by compressing the Hardy shift $M_z$ on $H^2_{\alpha,\beta}$. Formally, if $P_{H^2_{\alpha,\beta}}$ denote the orthogonal projection of $H^2(\D)$ onto $H^2_{\alpha, \beta}$, we  consider the operator $P_{H^2_{\alpha,\beta}}M_z|_{H^2_{\alpha, \beta}}$. We denote it by $S_{\alpha,\beta}$. The space $H^2_{\alpha, \beta}$ can also be regarded as a tridiagonal reproducing kernel Hilbert space with the set $\{\alpha+\beta z, z^n\quad n\geq 2\}$ as an orthonormal basis. Recall that, a functional Hilbert space is called a tridiagonal kernel space if it has an orthonormal basis of the form $\{f_n(z)=(a_n+b_nz)z^n\quad n\geq 0\}$, for some complex sequences $\{a_n\}_{n\geq 0}, \{b_n\}_{n\geq 0}$ (\cite{Adam 2001}). Here for $H^2_{\alpha, \beta}$, we have $a_0 =\alpha, b_0=\beta$, $a_1=b_n=0$, for $n\geq 1$ and $a_n=1$ for $n\geq 2$. However, we will not go in this direction.

In the next section, we show that, $S_{\alpha,\beta}$ can be viewed as a rank one perturbation of the Hardy shift on $H^2(\D)$ and at the same time, can be represented as the shift $M_z$ on some reproducing kernel Hilbert space (rkhs for short) of analytic functions (see Theorem \ref{shift}). This allows us to study the \textit{shifts} $S_{\alpha,\beta}$ in the light of perturbation theory and motivates the questions that we answer here. The analysis of perturbations of the form $M_z+F$ where $M_z$ is a shift and $F$ is some finite rank (or compact, Hilbert-Schmidt, etc.) operator has been studied extensively in the literature (e.g. \cite{ALP}, \cite{Clark}, \cite{Fuhrmann}, \cite{Nakamura}). The structure of their invariant subspaces and many essential properties like Normality, Hyponormality, etc. have received enormous attention (\cite{Curto 1}, \cite{Curto 2}, \cite{Curto}, \cite{DB}). While Substantial amount of work has been done on characterization of invariant subspaces of  perturbations, many of them~(\cite{ALP}, \cite{Nakamura}, \cite{Stampfli 1}, \cite{Tcaciuc}) evidently show, even dealing with rank one perturbations of shifts can be challenging. To study the invariant subspaces of a perturbation of form $S+F$ corresponding to the Hardy shift $S$, it is natural to look into the techniques and theories built for $S$. It is also natural to look for the connection between the corresponding structure and properties of $S+F$ and $S$. In our study, we establish a relation between the $S_{\alpha,\beta}$-invariant and Hardy shift-invariant subspaces. In fact, we show in section \ref{sec: Cyclicity}: \textsf{ For every nonzero, closed $S_{\alpha,\beta}$-invariant subspace $\clm$, there exists an inner function $\theta$ (unique upto a scalar multiple) such that either $\clm= z^2\theta H^2(\D)$ or $\dim (\clm\ominus z^2\theta H^2(\D))=1$} (see Theorem \ref{dim 1}). This is a key result in proving the cyclicity of $S_{\alpha,\beta}$-invariant subspaces. Based on this, we provide a representation of closed $S_{\alpha,\beta}$-invariant subspaces and show: \textsf{ Every $S_{\alpha,\beta}$-invariant subspace is cyclic} (see Theorem \ref{Isp char}). The rest of the paper is organized as follows:

\NI $\bullet$ In  Section \ref{sec: prep res}, we give some preparatory results that will be used throughout the paper. We also provide a necessary and sufficient condition on when two such $S_{\alpha,\beta}$ for different pairs of $\{\alpha,\beta\}$ are unitarily equivalent. The result is in terms of the simple equalities involving $\alpha,\beta$, which shows that the class of such non-equivalent $S_{\alpha,\beta}$ is quite large (see Proposition~\ref{Unitary equivalence S_ab}). We also find that the operators $S_{\alpha,\beta}$ in our study are hyponormal (\cite{Martin-Putinar}, \cite{Stampfli 2}). Following the results on cyclicity (section \ref{sec: Cyclicity}) as mentioned above,

\NI $\bullet$ Section \ref{sec: Characterization} deals with the characterization of $S_{\alpha,\beta}$-invariant subspaces. We prove a Beurling type analogue: \textsf{ For $\clm\subseteq H^2_{\alpha,\beta}$ being a nonzero, closed $S_{\alpha,\beta}$-invariant subspace, $\dim(\clm\ominus S_{\alpha,\beta}\clm)=1$.} (see Theorem \ref{codim 1}). Based on that, we provide a complete characterization of closed $S_{\alpha,\beta}$-invariant subspaces $\clm$ of $H^2_{\alpha,\beta}$ (Theorem \ref{S1 inv 1}).  Also, Theorem \ref{codim 1} provides a representation of a generating function of the one-dimensional subspace $(\clm\ominus S_{\alpha,\beta}\clm)$, which plays a crucial role in dealing with wandering subspace property (w.s.p.) in the next section.

\NI $\bullet$ In Section \ref{sec: Wandering Subspace property}, we discuss on the wandering subspaces of $S_{\alpha,\beta}$. We show, in contrast with the Hardy shift, $\ker S_{\alpha,\beta}^*$ need not generate the full space $H^2_{\alpha,\beta}$. Indeed, \textsf{ $\ker S_{\alpha,\beta}^*$ will generate $H^2_{\alpha,\beta}$ if and only if $|\alpha|^2\leq \frac{1}{u+1}$, $u$ being the unique real root of the cubic $z^3+3z^2+2z-1$.
} (see Theorem \ref{alpha less zero}). We provide supporting examples and counter examples in this regard.

In the later part of the same section, we consider the w.s.p. of the non-trivial closed $S_{\alpha,\beta}$-invariant subspaces $\clm$. We provide a necessary and sufficient condition for $\clm=[\clm\ominus S_{\alpha,\beta}\clm]_{S_{\alpha,\beta}}$ to hold (see Theorem \ref{r geq 1}) and clarify it with concrete examples. Eventually, this directs us to investigate the w.s.p. of $S_{\alpha,\beta}$- invariant subspaces with $|\beta|^2<\frac{1}{4+\gamma}$, where $\gamma$ is the unique real root of the cubic polynomial $z^3+7z^2+12z-1$. Along this line, we show that (unlike the case of classical shifts), there exist non-trivial closed $S_{\alpha,\beta}$-invariant subspaces $\clm\subseteq H^2_{\alpha,\beta}$ for which $\clm\neq [\clm\ominus S_{\alpha,\beta}\clm]$ ( see Theorem \ref{not wsp}).

\newsection{Preparatory Results}\label{sec: prep res}
In this section, we will provide some elementary results that will be used throughout. It is well-known (\cite{Paulsen Raghupati}): for $0<|\alpha|,|\beta|<1$ and $|\alpha|^2+|\beta|^2=1$,\ \ $H^2_{\alpha,\beta}$ is a closed subspace of codimension one of the Hardy space $H^2(\D)$ with $\{\alpha+\beta z, z^2, z^3, \ldots\}$ as an orthonormal basis.

Note that, $(\bar{\beta}-\bar{\alpha}z)\in H^2(\D)$ and
$\la \bar{\beta}-\bar{\alpha}z, \alpha+\beta z\ra= \overline{\beta\alpha}-\overline{\alpha\beta}=0$, $\la \bar{\beta}-\bar{\alpha}z, z^n\ra=0$ for all $n\geq 2$.
Again,
\begin{eqnarray}
1 &=& \beta(\bar{\beta}-\bar{\alpha}z)+\bar{\alpha}(\alpha+\beta z), \\ \label{1}
z &=& -\alpha(\bar{\beta}-\bar{\alpha}z)+\bar{\beta}(\alpha+\beta z).\label{z}
\end{eqnarray}
This implies,
$$\overline{span}\{\bar{\beta}-\bar{\alpha}z, \alpha+\beta z, z^2, z^3,\ldots\}= H^2(\D),\text{ and}$$
$$H^2(\D)\ominus H^2_{\alpha,\beta}=\C (\bar{\beta}-\bar{\alpha}z).$$
Hence, $\{\bar{\beta}-\bar{\alpha}z, \alpha+\beta z, z^2, z^3,\ldots\}$ forms an orthonormal basis of $H^2(\D)$.
\textsf{Let us denote $(\alpha+\beta z)$ by $f_0$.} We now proceed for the matrix representation of $S_{\alpha,\beta}= P_{H^2_{\alpha,\beta}}M_z|_{H^2_{\alpha,\beta}}$. Note that,
$$S_{\alpha,\beta}(f_0)= P_{H^2_{\alpha,\beta}}M_z(\alpha+\beta z)= P_{H^2_{\alpha,\beta}}(\alpha z+\beta z^2)= \alpha P_{H^2_{\alpha,\beta}}(z)+\beta z^2. $$
Since $(\bar{\beta}-\bar{\alpha}z)\perp H^2_{\alpha,\beta}$, it follows by \eqref{z}
\begin{equation}\label{S1_f0}
S_{\alpha,\beta}(f_0)= \alpha\bar{\beta}f_0+\beta z^2.
\end{equation}
Again for $n\geq 2$,
\begin{equation}\label{S1_z2}
S_{\alpha,\beta}(z^n)=P_{H^2_{\alpha,\beta}}M_z(z^n)=P_{H^2_{\alpha,\beta}}(z^{n+1})=z^{n+1}.
\end{equation}
Hence, with respect to the orthonormal basis $\{f_0,  z^n; n\geq 2\}$, $S_{\alpha,\beta}$ can be represented as
\begin{equation}\label{eqn: S_ab matrix}
[S_{\alpha,\beta}] = \begin{bmatrix}
\alpha\bar{\beta} & 0 & 0 & 0 & \cdots
\\
\beta & 0 & 0 & 0 & \cdots
\\
0 & 1 & 0 & 0 & \cdots
\\
0 & 0 & 1 & 0 & \cdots
\\
0 & 0 & 0 & 1 & \cdots
\\
\vdots & \vdots & \vdots& \vdots &\ddots
\end{bmatrix}.
\end{equation}
Again, with respect to the same orthonormal basis $\{f_0, z^2; n\geq 2\}$, $S_{\alpha,\beta}^*$ can be represented as
\begin{equation}\label{eqn: S_ab adjoint matrix}
[S_{\alpha,\beta}^*] = \begin{bmatrix}
\bar{\alpha}\beta & \bar{\beta} & 0 & 0 & 0 & \cdots
\\
0 & 0 & 1 & 0 & 0 & \cdots
\\
0 & 0 & 0 & 1 & 0 & \cdots
\\
0 & 0 & 0 & 0 & 1 & \cdots
\\
\vdots & \vdots & \vdots & \vdots &\vdots & \ddots
\end{bmatrix}.
\end{equation}

We now find criteria for unitary equivalence of the operators $S_{\alpha,\beta}$ on the spaces $H^2_{\alpha,\beta}$ for different pairs of $\{\alpha, \beta\}$. For that, we need the following Lemma on kernel of $S_{\alpha,\beta}^*$:
\begin{Lemma}\label{kernel S_ab adjoint}
For $\alpha, \beta$, both being nonzero and $|\alpha|^2+|\beta|^2=1$, $\ker S_{\alpha,\beta}^*= \mathbb{C} (f_0-\frac{\bar{\alpha}\beta}{\bar{\beta}}z^2)$.
\end{Lemma}
\begin{proof}
Clearly, we have by \eqref{eqn: S_ab adjoint matrix} $$S_{\alpha,\beta}^*(f_0-\frac{\bar{\alpha}\beta}{\bar{\beta}}z^2)=\bar{\alpha}\beta f_0-\bar{\alpha}\beta f_0=0.$$
Hence $(f_0-~\frac{\bar{\alpha}\beta}{\bar{\beta}}z^2)\in\ker S_{\alpha,\beta}^*$.
Again by \eqref{eqn: S_ab matrix}, $\text{range } S_{\alpha,\beta}=\text{span} \{\alpha\bar{\beta}f_0+\beta z^2, z^3, z^4,\ldots\}$. We show $\{f_0, z^2\}\subseteq \text{span}\{f_0-\frac{\bar{\alpha}\beta}{\bar{\beta}}z^2,\alpha\bar{\beta}f_0+\beta z^2\}$. Indeed,
$$z^2=\frac{1}{\beta(|\alpha|^2+1)}(\beta z^2+|\alpha|^2\beta z^2)=\frac{1}{\beta(|\alpha|^2+1)}\Big((\alpha\bar{\beta}f_0+\beta z^2)-\alpha\bar{\beta} (f_0-\frac{\bar{\alpha}\beta}{\bar{\beta}}z^2)\Big) \text{ and} $$
$f_0=(f_0-\frac{\bar{\alpha}\beta}{\bar{\beta}}z^2) + \frac{\bar{\alpha}\beta}{\bar{\beta}}z^2.$ Hence $(\text{ range }S_{\alpha,\beta})^{\perp}= \C(f_0-\frac{\bar{\alpha}\beta}{\bar{\beta}}z^2)=\ker S_{\alpha,\beta}^*$.
\end{proof}
For the pairs $(\alpha,\beta)$ and $(\alpha_1, \beta_1)$ on the unit sphere, it is well-known (\cite{Paulsen Raghupati}) that, $H^2_{\alpha,\beta}=H^2_{\alpha_1,\beta_1}$ if and only if $(\alpha,\beta)$ and $(\alpha_1, \beta_1)$ are scalar multiples of each other. Here we discuss the unitary equivalence of the operators $S_{\alpha,\beta}$ on $H^2_{\alpha,\beta}$ for different pairs of $\{\alpha, \beta\}$. Assume all $\alpha,\beta$ are nonzero and $(\alpha, \beta)\neq(\alpha_1, \beta_1)$. Then, we have the following proposition:
\begin{Proposition}\label{Unitary equivalence S_ab}
$S_{\alpha,\beta}$ on $H^2_{\alpha, \beta}$ and $S_{\alpha_1,\beta_1}$ on $H^2_{\alpha_1, \beta_1}$ are unitarily equivalent if and only if $|\alpha|=|\alpha_1|$; $|\beta|=|\beta_1|$ and $\frac{\alpha}{\alpha_1}=\frac{\bar{\beta_1}}{\bar{\beta}}$.
\end{Proposition}
\begin{proof}
Let us denote: $\tilde{f_0}=\alpha_1+\beta_1z$. As we noted earlier, the sets $B:=\{f_0, z^n; n\geq2\}$ and $B':=~\{\tilde{f_0}, z^n; n\geq2\}$ form orthonormal bases for $H^2_{\alpha, \beta}$ and $H^2_{\alpha_1, \beta_1}$ respectively.
Let there exist a unitary map $U: H^2_{\alpha,\beta}\longrightarrow H^2_{\alpha_1, \beta_1}$ such that
\begin{equation}\label{U begin}
US_{\alpha,\beta}=S_{\alpha_1,\beta_1}U
\end{equation}
We prove the necessary part in the following steps.

\NI\textsf{Step 1:} In this step, we show $U f_0= t_1\tilde{f_0}$ for some $|t_1|=1$. Note that, the matrix representations of $S_{\alpha_1,\beta_1}$ and $S_{\alpha_1,\beta_1}^*$ with respect to the orthonormal basis $\{\tilde{f_0}, z^n\}$ will be exactly same as that of $S_{\alpha,\beta}$ in (\eqref{eqn: S_ab matrix},\eqref{eqn: S_ab adjoint matrix}) with $\alpha_1, \beta_1$ in place of $\alpha, \beta$. By Lemma \ref{kernel S_ab adjoint}, $\ker S_{\alpha,\beta}^*=\C(f_0-\frac{\bar{\alpha}\beta}{\bar{\beta}}z^2)$ and hence $\ker S_{\alpha_1,\beta_1}^*=\C(\tilde{f_0}-\frac{\bar{\alpha_1}\beta_1}{\bar{\beta_1}}z^2)$.
Since, $S_{\alpha,\beta}^*U^*=U^*S_{\alpha_1,\beta_1}^*$ (By \eqref{U begin}), $U$ takes $\ker S_{\alpha,\beta}^*$ onto $\ker S_{\alpha_1,\beta_1}^*$. Hence, there exists $\gamma(\neq 0)\in\C$ such that \begin{equation}\label{U* ker to ker}
U^*(\tilde{f_0}-\frac{\bar{\alpha_1}\beta_1}{\bar{\beta_1}}z^2)=\gamma(f_0-\frac{\bar{\alpha}\beta}{\bar{\beta}}z^2).
\end{equation}
Note that, $I-S_{\alpha_1,\beta_1}^*S_{\alpha_1,\beta_1}= U (I- S_{\alpha,\beta}^*S_{\alpha,\beta})U^*$. By \eqref{eqn: S_ab matrix} and \eqref{eqn: S_ab adjoint matrix}, the matrix representations of $I-S_{\alpha,\beta}^*S_{\alpha,\beta}$ and $I-S_{\alpha,\beta}S_{\alpha,\beta}^*$ with respect to the orthonormal basis $\{f_0, z^n; n\geq2\}$ are

\begin{equation}\label{eqn: S_ab S_ab adjoint matrix}
[I-S_{\alpha,\beta}^*S_{\alpha,\beta}] = \begin{bmatrix}
|\alpha|^4 & 0 & 0 & 0 & \cdots
\\
0 & 0 & 0 & 0 & \cdots
\\
0 & 0 & 0 & 0 & \cdots
\\
0 & 0 & 0 & 0 & \cdots
\\
0 & 0 & 0 & 0 & \cdots
\\
\vdots & \vdots & \vdots& \vdots &\ddots
\end{bmatrix},\ \
[I-S_{\alpha,\beta}S_{\alpha,\beta}^*] = \begin{bmatrix}
1-|\alpha\beta|^2 & -\alpha\bar{\beta}^2 & 0 & 0 & \cdots
\\
-\bar{\alpha}\beta^2 & |\alpha|^2 & 0 & 0 & \cdots
\\
0 & 0 & 0 & 0 & \cdots
\\
0 & 0 & 0 & 0 & \cdots
\\
0 & 0 & 0 & 0 & \cdots
\\
\vdots & \vdots & \vdots& \vdots &\ddots
\end{bmatrix}.
\end{equation}
We continue to see the matrix representations of $(I-S_{\alpha_1,\beta_1}^*S_{\alpha_1,\beta_1})$ and $(I-S_{\alpha_1,\beta_1}S_{\alpha_1,\beta_1}^*)$ with respect to the basis $B'$ being same as \eqref{eqn: S_ab S_ab adjoint matrix} with entries $\alpha_1, \beta_1$ (replacing $\alpha, \beta$).
Let $U^*(\tilde{f_0})= r_0f_0+r_1z^2+r_2z^3+\cdots$, where $r_i\in\C$ for all $i$ and $\sum_{i=0}^{\infty}|r_i|^2=1$.
Now $(I-S_{\alpha_1,\beta_1}^*S_{\alpha_1,\beta_1})\tilde{f_0} = U(I-S_{\alpha,\beta}^*S_{\alpha,\beta})U^*\tilde{f_0}$ yields,
\begin{equation*}
|\alpha_1|^4\tilde{f_0}\ =\ U(I-S_{\alpha,\beta}^*S_{\alpha,\beta})(r_0f_0+r_1z^2+r_2z^3+\cdots)\ =\ U(r_0|\alpha|^4f_0)\quad (\text{ by }\eqref{eqn: S_ab S_ab adjoint matrix}).
\end{equation*}
Set $t_1= \frac{1}{r_0}|\frac{\alpha_1}{\alpha}|^4$. Then $|t_1|=1$ as $U$ is unitary and we have
\begin{equation}\label{Uf=t0}
U f_0= t_1\tilde{f_0}.
\end{equation}
\textsf{Step 2:} In this step, we show $Uz^2= t_2 z^2$ for some $|t_2|=1$.

\NI The equations \eqref{U* ker to ker} and \eqref{Uf=t0} together imply,
\begin{equation}\label{Uz^2}
U(f_0-\frac{\bar{\alpha}\beta}{\bar{\beta}}z^2)= t_1\tilde{f_0}-\frac{\bar{\alpha}\beta}{\bar{\beta}}U(z^2)=\frac{1}{\gamma}(\tilde{f_0}-\frac{\bar{\alpha_1}\beta_1}{\bar{\beta_1}}z^2).
\end{equation}
Since $\la Uz^2, \tilde{f_0}\ra= \la Uz^2, \bar{t_1}Uf_0 \ra= t_1 \la z^2, f_0\ra=0$, by \eqref{Uz^2} there exists $t_2\in\C$ such that
\begin{equation}\label{Ut2}
Uz^2= t_2 z^2,\ \ \text{ where } |t_2|=1.
\end{equation}

\NI\textsf{Step 3:} In this step, we aim to show $U(z^3)=t_3 z^3$ for some scalar $t_3$ with $|t_3|=1$.

\NI A simple computation will reveal that the matrix representations of $(I-S_{\alpha,\beta}^2S_{\alpha,\beta}^{*2})$ and $(I-S_{\alpha_1,\beta_1}^2S_{\alpha_1,\beta_1}^{*2})$ with respect to the orthonormal bases $B$ and $B'$ respectively are given by,
\begin{equation}\label{eqn: S_ab square}
[I-S_{\alpha,\beta}^2S_{\alpha,\beta}^{*2}] = \begin{bmatrix}

1-|\alpha\bar{\beta}|^4 & -\bar{\alpha}|\beta|^2(\alpha\bar{\beta})^2 & -\bar{\beta}(\alpha\bar{\beta})^2 & 0 & \dots
\\
-\alpha|\beta|^2(\bar{\alpha}\beta)^2 & 1-|\alpha|^2|\beta|^4 & -\alpha|\beta|^2\bar{\beta} & 0 & \ddots
\\
-\beta(\bar{\alpha}\beta)^2 & -\bar{\alpha}\beta|\beta|^2 & 1-|\beta|^2 & 0 & \ddots
\\
0 & 0 & 0 & 0 & \ddots
\\
\vdots & \vdots & \vdots&\vdots &\ddots
\end{bmatrix}
\end{equation}

\begin{equation}\label{eqn: S_ab prime square}
[I-S_{\alpha_1,\beta_1}^2S_{\alpha_1,\beta_1}^{*2}] = \begin{bmatrix}

1-|\alpha_1\bar{\beta_1}|^4 & -\bar{\alpha_1}|\beta_1|^2(\alpha_1\bar{\beta_1})^2 & -\bar{\beta_1}(\alpha_1\bar{\beta_1})^2 & 0 & \dots
\\
-\alpha_1|\beta_1|^2(\bar{\alpha_1}\beta_1)^2 & 1-|\alpha_1|^2|\beta_1|^4 & -\alpha_1|\beta_1|^2\bar{\beta_1} & 0 & \ddots
\\
-\beta_1(\bar{\alpha_1}\beta_1)^2 & -\bar{\alpha_1}\beta_1|\beta_1|^2 & 1-|\beta_1|^2 & 0 & \ddots
\\
0 & 0 & 0 & 0 & \ddots
\\
\vdots & \vdots & \vdots&\vdots &\ddots
\end{bmatrix}.
\end{equation}
Again, since $(I-S_{\alpha_1,\beta_1}^2S_{\alpha_1,\beta_1}^{*2})= U(I-S_{\alpha,\beta}^2S_{\alpha,\beta}^{*2})U^*$ and $U^*z^3$ is orthogonal to $f_0$ and $z^2$, there exist scalars $r_i$\ ($i\geq 0$) with $\sum_{i\geq 0}|r_i|^2=1$ such that $U^*z^3= r_0z^3+r_1z^4+r_2z^5+\cdots$. Hence
\[
\begin{split}
(I-S_{\alpha_1,\beta_1}^2S_{\alpha_1,\beta_1}^{*2})z^3 &=U(I-S_{\alpha,\beta}^2S_{\alpha,\beta}^{*2})U^*(z^3) \\
&= U(I-S_{\alpha,\beta}^2S_{\alpha,\beta}^{*2})(r_0z^3+r_1z^4+r_2z^5+\cdots)\\
&=r_0 U(-\bar{\beta}(\alpha\bar{\beta})^2 f_0-\alpha|\beta|^2\bar{\beta}z^2+ (1-|\beta|^2)z^3)\ \ (\text{ by \eqref{eqn: S_ab square}}) \\
&= r_0 \Big((-\bar{\beta}(\alpha\bar{\beta})^2)t_1 \tilde{f_0}- \alpha|\beta|^2\bar{\beta}t_2 z^2+ (1-|\beta|^2)U(z^3)\Big) (\text{ by \eqref{Uf=t0}, \eqref{Ut2}}).
\end{split}
\]
Now, $(I-S_{\alpha_1,\beta_1}^2S_{\alpha_1,\beta_1}^{*2})z^3= -\bar{\beta_1}(\alpha_1\bar{\beta_1})^2)\tilde{f_0}- \alpha_1|\beta_1|^2\bar{\beta_1} z^2+ (1-|\beta_1|^2)z^3$ (by \eqref{eqn: S_ab prime square}) and again, as $U(z^3)$ is orthogonal to $\tilde{f_0}(=\bar{t_1}U(f_0))$ and $z^2(=\bar{t_2} U(z^2))$, by above computation it follows that, there exists $t_3\in\C$ with $|t_3|=1$ such that
\begin{equation}\label{U z3=t3}
Uz^3=t_3z^3.
\end{equation}
\NI\textsf{Step 4:} In this step we prove the following claim:

\NI\textsf{ Claim:} For $n\geq 1$,
\begin{eqnarray}
S_{\alpha,\beta}^n f_0 &=& (\alpha\bar{\beta})^nf_0+\beta\big[(\alpha\bar{\beta})^{n-1}z^2+(\alpha\bar{\beta})^{n-2}z^3+\cdots+(\alpha\bar{\beta})z^n+ z^{n+1}\big]\ \ \text{and} \label{S_ab power n}
\\
S_{\alpha,\beta} ^n z^m &=& z^{m+n},\quad  m\geq 2.\label{S1 zm}
\end{eqnarray}

Clearly, $S_{\alpha,\beta}z^m= z^{m+1}$ for $m\geq 2$ (by \eqref{eqn: S_ab matrix}) and hence for $n\geq 1$, by successive application of $S_{\alpha, \beta}$, we will have $S_{\alpha,\beta} ^n z^m = z^{m+n}$ for all $ m\geq 2$. This settles \eqref{S1 zm}. We prove \eqref{S_ab power n} by induction.
Note by \eqref{eqn: S_ab matrix}, $S_{\alpha,\beta} f_0 = \alpha\bar{\beta}f_0+\beta z^2$ and hence
\[
S_{\alpha,\beta}^2f_0 = S_{\alpha,\beta}(\alpha\bar{\beta}f_0+\beta z^2)=\alpha\bar{\beta}(\alpha\bar{\beta}f_0+\beta z^2)+\beta z^3 =(\alpha\bar{\beta})^2f_0+\beta(\alpha\bar{\beta}z^2+z^3),
\]
showing \eqref{S_ab power n} holds for $n = 1,2$. Now for the induction step, let it hold for some $n\geq 3$. Then the claim \eqref{S_ab power n} is settled by the following computation:
\[
\begin{split}
S_{\alpha,\beta}^{n+1}f_0 &= S_{\alpha,\beta} S_{\alpha,\beta}^nf_0  \\
&= S_{\alpha,\beta} \Big[(\alpha\bar{\beta})^nf_0+\beta\big((\alpha\bar{\beta})^{n-1}z^2+(\alpha\bar{\beta})^{n-2}z^3+\cdots+(\alpha\bar{\beta})z^n+ z^{n+1}\big)\Big]\\
&=(\alpha\bar{\beta})^n(\alpha\bar{\beta}f_0+\beta z^2)+\beta\big[(\alpha\bar{\beta})^{n-1}z^3+(\alpha\bar{\beta})^{n-2}z^4+\cdots+(\alpha\bar{\beta})z^{n+1}+ z^{n+2}\big] \\
&=(\alpha\bar{\beta})^{n+1}f_0 + \beta\big[(\alpha\bar{\beta})^n z^2+(\alpha\bar{\beta})^{n-1}z^3+(\alpha\bar{\beta})^{n-2}z^4+\cdots+(\alpha\bar{\beta})z^{n+1}+ z^{n+2}\big].
\end{split}
\]

\NI\textsf{Step 5:} In this step, we show by induction, there exist scalars $t_n\in\C$ with $|t_n|=1$ such that
\begin{equation}\label{U diagonal}
U(z^n)=t_nz^n\quad \text{for all } n\geq 2.
\end{equation}
\NI The cases $n=2, 3$ follows by \eqref{Ut2} and \eqref{U z3=t3}. Assume that, $U(z^n)=t_n z^n$ with $|t_n|=1$ hold for $n= 2, 3, \ldots, m$ for some $m \geq 4$. We show, $U(z^{m+1})=t_{m+1}z^{m+1}$ for some $|t_{m+1}|=1$.
Note that, for all $m\geq 1$
\begin{equation}\label{S_ab m power matrix}
(I-S_{\alpha_1,\beta_1}^mS_{\alpha_1,\beta_1}^{*m})=U(I-S_{\alpha,\beta}^mS_{\alpha,\beta}^{*m})U^*.
\end{equation}
We now consider matrix representation of $(I- S_{\alpha,\beta}^m S_{\alpha,\beta}^{*m})$, for $m\geq 3$ with respect to the orthonormal basis $B:=\{f_0; z^n, n\geq2\}$. It follows by \eqref{S_ab power n}, there exist complex scalars $a_{i,j}, 1\leq i,j\leq (m+1)$ with $a_{1,1}=1-|\alpha\bar{\beta}|^{2m}$ and $a_{m+1, m+1}=1-|\beta|^2$ such that
\begin{equation}\label{eqn: S power m}
[I-S_{\alpha,\beta}^m S_{\alpha,\beta}^{*m}] = \begin{bmatrix}
1-|\alpha\bar{\beta}|^{2m} & a_{12} & \dots & a_{1,m+1} & 0 &  \dots
\\
a_{21} & a_{22} & \dots & a_{2, m+1} & 0 &  \dots
\\
\vdots & \vdots & \dots & \vdots & \vdots & \dots
\\
a_{m+1,1} & a_{m+1, 2} & \dots  & 1-|\beta|^2 & 0& \dots
\\
0 & 0 & \dots & \dots  & 0& \dots
\\
\vdots & \vdots & \vdots &\vdots &\vdots &\dots
\end{bmatrix}.
\end{equation}
The matrix representation of $(I-S_{\alpha_1,\beta_1}^mS_{\alpha_1,\beta_1}^{*m})$ for $m\geq 3$ with respect to the basis $B'=\{\tilde{f_0}; z^n, n\geq2\}$, will also have similar form:
\begin{equation}\label{eqn: S prime power m}
[I-S_{\alpha_1,\beta_1}^m S_{\alpha_1,\beta_1}^{*m}] = \begin{bmatrix}
1-|\alpha_1\bar{\beta_1}|^{2m} & a_{12}' & \dots & a_{1,m+1}' & 0 &  \dots
\\
a_{21}' & a_{22}' & \dots & a_{2, m+1}' & 0 &  \dots
\\
\vdots & \vdots & \dots & \vdots & \vdots & \dots
\\
a_{m+1,1}' & a_{m+1, 2}' & \dots  & 1-|\beta_1|^2 & 0 & \dots
\\
0 & 0 & \dots & \dots  & 0& \dots
\\
\vdots & \vdots & \vdots &\vdots &\vdots &\dots
\end{bmatrix},
\end{equation}
for some scalars $a_{i,j}', 1\leq i,j\leq (m+1)$. Since $U^*(z^{m+1})$ is orthogonal to $f_0, z^2,\ldots z^m$, there exist scalars $r_i$, ($i=0,1, 2, \ldots$) with $\sum_{i\geq 0}|r_i|^2=1$ such that
$U^*z^{m+1}=(r_0z^{m+1}+r_1z^{m+2}+~\cdots)$. Hence, by \eqref{S_ab m power matrix}
\[
\begin{split}
(I-S_{\alpha_1,\beta_1}^mS_{\alpha_1,\beta_1}^{*m})z^{m+1} &= U(I-S_{\alpha,\beta}^mS_{\alpha,\beta}^{*m})(r_0z^{m+1}+r_1z^{m+2}+\cdots)  \\
&= r_0 U \Big(a_{1,m+1}f_0+a_{2, m+1}z^2+\cdots +(1-|\beta|^2)z^{m+1}\Big)\ (\text{ by \eqref{eqn: S power m}}) \\
&=r_0(a_{1,m+1}t_1\tilde{f_0}+a_{2, m+1}t_2z^2+\cdots +(1-|\beta|^2)Uz^{m+1}),
\end{split}
\]
where the last equality follows by Step $1$ and induction hypothesis. Since
$$(I-S_{\alpha_1,\beta_1}^mS_{\alpha_1,\beta_1}^{*m})z^{m+1}= (a_{1,m+1}'\tilde{f_0}+a_{2, m+1}'z^2+\cdots +(1-|\beta_1|^2)z^{m+1}),$$
by similar argument as in Step $2$ and Step $3$, we must have $U(z^{m+1})=t_{m+1}z^{m+1}$ for some scalar $t_{m+1}$ satisfying $|t_{m+1}|=1$.

\NI\textsf{Step 6:} In this step, we consider the matrix representation of $U$.

\NI Note that, by \eqref{Uf=t0} and \eqref{U diagonal} and with respect to the orthonormal bases $B$ on $H^2_{\alpha,\beta}$ and $B'$ on $H^2_{\alpha_1,\beta_1}$, $U$ can be represented as
\begin{equation}\label{U matrix}
\begin{bmatrix}
t_1& 0 & 0 & 0 & \dots
\\
0 & t_2 & 0 & 0 & \ddots
\\
0 & 0 & t_3& 0 & \ddots
\\
0 & 0 & 0 & t_4 & \ddots
\\
\vdots & \vdots & \vdots&\ddots &\ddots
\end{bmatrix},
\end{equation}
where $|t_n|=1$ for all $n\geq 1$.

\NI\textsf{Final Step :} We now show $|\alpha_1|=|\alpha|$, $|\beta_1|=|\beta|\text{ and } \frac{\alpha}{\alpha_1}=\frac{\bar{\beta_1}}{\bar{\beta}}$.

\NI The matrix representations \eqref{eqn: S_ab matrix}, \eqref{U matrix} together with the relation $S_{\alpha_1,\beta_1}=US_{\alpha,\beta}U^*$ yield,
\begin{equation}\label{S_ab prime = S_ab}
\begin{bmatrix}
\alpha_1\bar{\beta_1}& 0 & 0 & 0 & \dots
\\
\beta_1 & 0 & 0 & 0 & \ddots
\\
0 & 1 & 0 & 0 & \ddots
\\
0 & 0 & 1 & 0 & \ddots
\\
0 & 0 & 0 & 1 & \ddots
\\
\vdots & \vdots & \vdots&\ddots &\ddots
\end{bmatrix}
= \begin{bmatrix}
\alpha\bar{\beta}& 0 & 0 & 0 & \dots
\\
\bar{t_1}t_2\beta & 0 & 0 & 0 & \ddots
\\
0 & \bar{t_2}t_3 & 0 & 0 & \ddots
\\
0 & 0 & \bar{t_3}t_4 & 0 & \ddots
\\
0 & 0 & 0 & \bar{t_4}t_5 & \ddots
\\
\vdots & \vdots & \vdots&\ddots &\ddots
\end{bmatrix}.
\end{equation}
Clearly by \eqref{S_ab prime = S_ab}, $\alpha_1\bar{\beta_1}=\alpha\bar{\beta}$,\ \ $\beta_1=\bar{t_1}t_2\beta$\ \ and\ \ $\bar{t_n}t_{n+1}=1$ for all $n\geq 2$. Since $|t_n|^2=t_n\bar{t_n}=1$ for all $n\geq 1$, we will have $t_n=t_{n+1}=t_2$\ \ for all $n\geq 2$. Therefore, the unitary operator $U$ can be defined as: $U(f_0)=t_1f_0$\ \ and\ \ $U(z^n)=t_2z^n$ for $n\geq 2$. Clearly, $|\beta_1|=|t_2\bar{t_1}\beta|=|\beta|$, which together with $\alpha_1\bar{\beta_1} =\alpha\bar{\beta}$ implies $|\alpha|=|\alpha_1|$ and $\frac{\alpha}{\alpha_1}=\frac{\bar{\beta_1}}{\bar{\beta}}$.

Conversely, let us assume $$|\alpha_1|=|\alpha|,\quad |\beta_1|=|\beta|\text{ and } \frac{\alpha}{\alpha_1}=\frac{\bar{\beta_1}}{\bar{\beta}}.$$
We define an operator $U: H^2_{\alpha, \beta}\longrightarrow H^2_{\alpha_1, \beta_1}$ by,
$$U(f_0)=\tilde{f_0}\ \ \text{ and }\ \ U(z^n)=\frac{\beta_1}{\beta} z^n,\quad  n\geq 2.$$
 Then $U$ is unitary (as $|\beta_1|=|\beta|$) and the proof is complete by the below computations:
\[
US_{\alpha,\beta}U^*(\tilde{f_0})=US_{\alpha,\beta}(f_0)=U(\alpha\bar{\beta}f_0+\beta z^2)=\alpha\bar{\beta}\tilde{f_0}+\beta\frac{\beta_1}{\beta}z^2=\alpha_1\bar{\beta_1}\tilde{f_0}+\beta_1z^2=S_{\alpha_1,\beta_1}(\tilde{f_0}),
\]
where the fourth equality follows by the assumption $\frac{\alpha}{\alpha_1}=\frac{\bar{\beta_1}}{\bar{\beta}}$,
and finally, for $n\geq 2$
\[
US_{\alpha,\beta}U^*(z^n)=US_{\alpha,\beta}(\frac{\beta}{\beta_1} z^n)= \frac{\beta}{\beta_1}U(z^{n+1})=\frac{\beta}{\beta_1}\frac{\beta_1}{\beta} z^{n+1}=S_{\alpha_1,\beta_1}(z^n).
\]
\end{proof}
\begin{Remark}
\begin{enumerate}
\item The proof of Proposition \ref{Unitary equivalence S_ab} shows that, if $U$ is a unitary map intertwining $S_{\alpha,\beta}$ for different pairs of $\{\alpha,\beta\}$, then $U$ is a diagonal
operator with at most two distinct diagonal entries i.e., $ U:=\diag\big(t_1, t, t,\cdots\big)$ for $t_1, t\in \C$ with $|t_1|=|t|=1$.

\item The operator $S_{\alpha,\beta}$ is hyponormal. In fact, one can write
$$S_{\alpha,\beta}^*S_{\alpha,\beta}-S_{\alpha,\beta}S_{\alpha,\beta}^*=(I-S_{\alpha,\beta}S_{\alpha,\beta}^*)-(I-S_{\alpha,\beta}^*S_{\alpha,\beta}),$$
and then by \eqref{eqn: S_ab S_ab adjoint matrix}, it follows that $S_{\alpha,\beta}^*S_{\alpha,\beta}-S_{\alpha,\beta}S_{\alpha,\beta}^*\geq0$.
\end{enumerate}
\end{Remark}
The Proposition \ref{Unitary equivalence S_ab} establishes the existence of a large class of non-equivalent $S_{\alpha,\beta}$.
We now show that any such $S_{\alpha,\beta}$ on $H^2_{\alpha,\beta}$ can actually be viewed as the multiplication by $z$ on some reproducing kernel Hilbert space. Recall that, a bounded linear operator $T$ is left invertible if $T^*T$ is invertible and analytic if $\bigcap_{n\geq 1}T^n\clh=\{0\}$. It is well-known (\cite{SS}): \textsf{If a bounded linear operator $T$ on a Hilbert space $\clh$ is left-invertible and analytic, then it is unitarily equivalent to the shift operator $M_z$ on some rkhs $\clh_k$.} In our context i.e., for nonzero $\alpha, \beta$ with $|\alpha|^2+|\beta|^2=1$, we show in below theorem, $S_{\alpha,\beta}$ is left invertible and analytic.
\begin{Theorem}\label{shift}
$S_{\alpha,\beta}$ on $H^2_{\alpha,\beta}$ is a shift $M_z$ on some rkhs $\clh_k$.
\end{Theorem}
\begin{proof}
By\eqref{eqn: S_ab matrix} and \eqref{eqn: S_ab adjoint matrix},  it follows that the matrix representation of $S_{\alpha,\beta}^*S_{\alpha,\beta}$ with respect to the orthonormal basis $\{f_0, z^n: n\geq 0\}$ is given by,
\begin{equation}\label{eqn: S_ab adjoint S_ab matrix}
[S_{\alpha,\beta}^*S_{\alpha,\beta}] = \begin{bmatrix}
|\alpha\beta|^2+|\beta|^2 & 0 & 0 & 0 & \cdots
\\
0 & 1 & 0 & 0 & \cdots
\\
0 & 0 & 1 & 0 & \cdots
\\
0 & 0 & 0 & 1 & \cdots
\\
\vdots & \vdots & \vdots& \vdots &\ddots
\end{bmatrix}.
\end{equation}
Since $|\alpha\beta|^2+|\beta|^2\neq 0$, $S_{\alpha,\beta}^*S_{\alpha,\beta}$ is invertible and hence, $S_{\alpha,\beta}$ is left invertible.

Let $U: H^2_{\alpha,\beta}\longrightarrow H^2(\D)$ be the canonical unitary map defined by $U(f_0)=1$ and $U(z^n)=z^{n-1}$ for all $n\geq 2$. Then it is easy to see via $U$, the operator $S_{\alpha,\beta}$ on $H^2_{\alpha,\beta}$ can be viewed as the rank one perturbation of the Hardy shift having the same matrix representation as \eqref{eqn: S_ab matrix} with respect to the standard orthonormal basis $\{z^n\}_{n\geq 0}$ on $H^2(\D)$. More precisely, $US_{\alpha,\beta}U^*= S+F$ where $S$ is the unilateral shift of multiplicity one acting on $\{z^n\}_{n\geq 0}$ and $F$ on $H^2(\D)$ is defined as $F(1)=\alpha\bar{\beta}+(\beta-1) z$ and $F(z^n)=0$ for all $n\geq 1$.
To show the analyticity of $S_{\alpha,\beta}$, we need the following result (Corollary $4.5.$ in \cite{Das}):
\textsf{ If $T$ is a bounded linear operator on $H^2(\D)$ with $T=S_k+F$ where~ $S_k$ acts on $\{z^n\}_{n\geq 0}$ as the unilateral shift of multiplicity $1\leq k <\infty$ and $F$ is defined by,
\[
F(z^n)=\begin{cases}
\alpha_0+\alpha_1 z+\cdots+(\alpha_k-1)z^k & \mbox{if } n=0 \\
0 & \mbox{if } n\geq 1,
\end{cases}
\]
for $\{\alpha_j\}_{j=0}^k\in\C$ and $0<|\alpha_0|\leq 1$, then $T$ is analytic if and only if $\alpha_j\neq 0$ for some $1\leq j\leq k$.}

\NI Since, in our setting $0<|\alpha\beta|<1$ and $\beta\neq 0$, it readily follows that, $S_{\alpha,\beta}$ is analytic. This completes the proof as the property of being analytic is a unitary invariant.
\end{proof}
\begin{Remark}
The analyticity part of $S_{\alpha,\beta}$ in Theorem \ref{shift} can also be proved by using a result from (\cite{SS}): \textsf{A left invertible operator $T\in\clb(\clh)$ is analytic if and only if it's Cauchy dual $T'=T(T^*T)^{-1}$ satisfies $[\ker T'^*]_{T'}=\clh$.} In our setting, $S_{\alpha,\beta}'= S_{\alpha,\beta}(S_{\alpha,\beta}^*S_{\alpha,\beta})^{-1}$ and $\ker S_{\alpha,\beta}'^*=\ker S_{\alpha,\beta}^*=\C(f_0-\frac{\bar{\alpha}\beta}{\bar{\beta}}z^2)$ (Lemma \ref{kernel S_ab adjoint}). A simple computation shows that, the matrix representation of $S_{\alpha,\beta}'$ with respect to the orthonormal basis $\{f_0, z^n, n\geq 2\}$ of $H^2_{\alpha,\beta}$ is
\begin{equation}\label{eqn: S_2 matrix}
[S_{\alpha,\beta}'] = \begin{bmatrix}
\frac{\alpha\bar{\beta}}{p} & 0 & 0 & 0 & \cdots
\\
\frac{\beta}{p} & 0 & 0 & 0 & \cdots
\\
0 & 1 & 0 & 0 & \cdots
\\
0 & 0 & 1 & 0 & \cdots
\\
0 & 0 & 0 & 1 & \cdots
\\
\vdots & \vdots & \vdots& \vdots &\ddots
\end{bmatrix},
\end{equation}
where $p=|\alpha\beta|^2+|\beta|^2$. Let us denote $f_2:=(f_0-\frac{\bar{\alpha}\beta}{\bar{\beta}}z^2)$. We show $[f_2]_{S_{\alpha,\beta}'}=H^2_{\alpha, \beta}$. Note that, $[f_2]_{S_{\alpha,\beta}'}=\C f_2 + [S_{\alpha,\beta}'f_2- \frac{\bar{\alpha}{\beta}}{p}f_2]_{S_{\alpha,\beta}'}$. By \eqref{eqn: S_2 matrix},
$S_{\alpha,\beta}'f_2- \frac{\bar{\alpha}{\beta}}{p}f_2= z^2\big(\frac{1}{\bar{\alpha}\beta}-z\big)$. Again, $|\alpha\beta|<1$ and $S_{\alpha,\beta}'(z^n)=z^{n+1}$ ($n\geq2$), together imply $[S_{\alpha,\beta}'f_2- \frac{\bar{\alpha}{\beta}}{p}f_2]_{S_{\alpha,\beta}'}= z^2H^2(\D)$ and hence $[f_2]_{S_{\alpha,\beta}'}= H^2_{\alpha,\beta}$. Consequently, $S_{\alpha,\beta}$ is analytic.
\end{Remark}

\newsection{ Cyclicity}\label{sec: Cyclicity}
We now turn to the invariant subspaces of $S_{\alpha,\beta}$ and the related properties. One of the most well-studied properties about the invariant subspaces of a shift operator is cyclicity. Needless to say, the cyclicity is a complex problem in shift-invariant subspaces and the corresponding perturbations (\cite{ALP, Sola et al}). In this section, we discuss on cyclicity of $S_{\alpha,\beta}$-invariant subspaces and prove that, any $S_{\alpha,\beta}$-invariant subspace is cyclic. As a first step, we prove the following relation between $S_{\alpha,\beta}$-invariant and Hardy shift $S$-invariant subspaces.
\begin{Lemma}\label{shift isp included}
Let $\clm$ be a nonzero, closed subspace of $H^2_{\alpha,\beta}$ such that $S_{\alpha,\beta}\clm\subseteq\clm$. Then there exists an inner function $\theta$ such that $z^2\theta H^2(\D)\subseteq\clm$.
\end{Lemma}
\begin{proof}
Since $\clm\neq \{0\}$, there exists $f\in\clm$ such that $f\neq 0$. Two cases can arise:

\NI\textbf{Case 1:} $f(0)=0$.

\NI Since $f\in H^2_{\alpha,\beta}$, we can write $f=c_0f_0+z^2g= c_0(\alpha+\beta z)+z^2g$, for some scalar $c_0\in\C$ and $g\in H^2(\D)$.
Now, $f(0)=0$ implies $c_0=0$ (as $\alpha\neq 0$) and hence $f=z^2g$. Since $f\neq 0$, $g$ must be nonzero. Now by \eqref{eqn: S_ab matrix}, $S_{\alpha,\beta}(f)=S_{\alpha,\beta}(z^2g)=z^3g=S(z^2g)$ and hence by Beurling's Theorem, there exists an inner function $\theta$ such that, $[f]_{S_{\alpha,\beta}}=[f]_S=z^2\theta H^2(\D)\subseteq\clm$.

\NI\textbf{Case 2:} $f(0)\neq 0$.

\NI We can write $f=c_0f_0+z^2g$ for some nonzero $c_0\in\C$ and $g\in H^2(\D)$. Then $$S_{\alpha,\beta}f= c_0(\alpha\bar{\beta}f_0+\beta z^2)+z^3g\ \ ( \text{ by } \eqref{eqn: S_ab matrix}),$$ and hence
\begin{equation}\label{S1 f non zero}
S_{\alpha,\beta}f-\alpha\bar{\beta}f = c_0(\alpha\bar{\beta}f_0+\beta z^2)+z^3g-\alpha\bar{\beta}(c_0f_0+z^2g) = z^2 \big(c_0\beta+(z-\alpha\bar{\beta})g\big).
\end{equation}
Note that, $c_0\beta+(z-\alpha\bar{\beta})g\neq 0$. Otherwise we will have $(z-\alpha\bar{\beta})g=-c_0\beta$, which would imply $(z-\alpha\bar{\beta})$ is invertible in $H^2(\D)$ with $|\alpha\beta|<1$. Now by \eqref{S1 f non zero} and \eqref{eqn: S_ab matrix}, for all $n\geq 1$
$$S_{\alpha,\beta}^n(S_{\alpha,\beta}f-\alpha\bar{\beta}f)=z^{n+2}(c_0\beta+(z-\alpha\bar{\beta})g)=S^n(S_{\alpha,\beta}f-\alpha\bar{\beta}f)\in [f]_{S_{\alpha,\beta}}\subseteq\clm.
$$
Hence, there exists an inner function $\theta$ such that $z^2\theta H^2(\D)\subseteq[f]_{S_{\alpha,\beta}}\subseteq \clm$.
\end{proof}
\begin{Corollary}\label{irreducible}
The operator $S_{\alpha,\beta}$ on $H^2_{\alpha, \beta}$ is irreducible.
\end{Corollary}
\begin{proof}
Suppose on the contrary, $\clm (\neq 0)\subseteq H^2_{\alpha,\beta}$ be a proper closed subspace such that $S_{\alpha,\beta}\clm\subseteq\clm$ and $S_{\alpha,\beta}\clm^{\perp}\subseteq\clm^{\perp}$. By Lemma \ref{shift isp included}, there exist inner functions $\theta_1, \theta_2$ such that $z^2\theta_1 H^2(\D)\subseteq\clm$, $z^2\theta_2 H^2(\D)\subseteq\clm^{\perp}$. Then $z^2\theta_1\theta_2\in \clm\bigcap\clm^{\perp}$, which is a contradiction.
\end{proof}
We are now in a position to state the following property of $S_{\alpha,\beta}$-invariant subspaces. This will play a crucial role in proving their cyclicity.
\begin{Theorem}\label{dim 1}
Let $\clm\subseteq H^2_{\alpha,\beta}$ be a nonzero closed $S_{\alpha,\beta}$-invariant subspace. Then there exists an inner function $\theta$ such that either $\clm=z^2\theta H^2(\D)$ or $\dim \big(\clm\ominus z^2\theta H^2(\D)\big)=1$. Moreover $\theta$ is unique upto a scalar multiple of unit modulus.
\end{Theorem}
\begin{proof}
Let
\begin{equation}\label{M1}
\clm_1= \overline{span}\{f\in \clm: f(z)= z^2g(z)\text{ for some } g\in H^2(\D)\}.
\end{equation}
Since $\clm\neq 0$, similar argument as in Lemma \ref{shift isp included} will show, $\clm_1\neq\{0\}$. Clearly, $\clm_1$ is a closed linear subspace of $\clm$ and every element has a zero of order $\geq2$ at $0$. Note that, for any $f\in\clm_1$,
$$S_{\alpha,\beta}f=S_{\alpha,\beta}(z^2g)=z^3g=S(z^2g)=S(f).$$
Hence, $\clm_1$ is a closed subspace of $H^2(\D)$ invariant under $S$ (and also under $S_{\alpha,\beta}$) and by Beurling's Theorem \cite{Beurling}, there exists an inner function $\theta$ such that
\begin{equation}\label{M1 first}
\clm_1=z^2\theta H^2(\D).
\end{equation}
Note that, if $\clm= H^2_{\alpha,\beta}$, then $\clm_1= z^2H^2(\D)$ and hence, $\clm\ominus z^2H^2(\D)=\C f_0$ and $\dim(\clm\ominus z^2H^2(\D))=1$. So, let us assume that $\clm$ is proper. We have the following two cases:

\NI \textbf{Case 1:} $f(0)=0$ for all $f\in\clm$.

\NI If $f(0)=0$ for all $f\in\clm$, we will have $\clm=\clm_1=z^2\theta H^2(\D)$, for some inner function $\theta$. Also, uniqueness of $\theta$ is clear in this case.

\NI \textbf{Case 2:} $f(0)\neq 0$ for some $f\in\clm$.

\NI Then by \eqref{M1 first}, $\clm\ominus\clm_1=\clm\ominus z^2\theta H^2(\D)\neq \{0\}$. Let $f\in \clm\ominus z^2\theta H^2(\D)$. We can take $f$ (upto a scalar multiple) as $f=f_0+z^2g$ for some $g\in H^2(\D)$. Note that, $g\neq 0$. Otherwise, $[f]_{S_{\alpha,\beta}}=[f_0]_{S_{\alpha,\beta}}\subseteq \clm$ and $[f_0]_{S_{\alpha,\beta}}=\overline{span}\{f_0, \alpha\bar{\beta}f_0+\beta z^2,\ldots\}= H^2_{\alpha,\beta}$ (by \eqref{eqn: S_ab matrix}) would imply $\clm= H^2_{\alpha,\beta}$, contradicting the assumption: $\clm$ is proper.

Now proceeding exactly as in \eqref{S1 f non zero} with $c_0=1$, we have
$$S_{\alpha,\beta}f-\alpha\bar{\beta}f= z^2(\beta+(z-\alpha\bar{\beta})g)\in\clm_1 (\text{ by } \eqref{M1}),$$
and hence by \eqref{M1 first}, there exists $h (\neq 0)\in H^2(\D)$ such that
$$z^2(\beta+(z-\alpha\bar{\beta})g)= z^2\theta h,$$ which is equivalent to
\begin{equation}\label{eqn: beta+g=theta h}
\beta+(z-\alpha\bar{\beta})g = \theta h\quad \big(\text{ as } \beta+(z-\alpha\bar{\beta})g\neq 0\big).
\end{equation}
Taking the inner product with $\theta$ on both sides of $\eqref{eqn: beta+g=theta h}$, we have
\begin{equation}\label{eqn: h(0)= theta+ zg}
\beta\overline{\theta(0)}+ \la zg, \theta\ra-\alpha\bar{\beta}\la g, \theta\ra =\la h, 1\ra.
\end{equation}
As $f(=f_0+z^2g) \perp z^2\theta H^2(\D)$, we have $\la f, z^{n+2}\theta\ra=0$ for all $n\geq 0$. This yields,
\begin{equation}\label{g perp theta Hardy}
\la g, z^n\theta\ra=0\text{ for all } n\geq 0.
\end{equation}
Hence, the equations $\eqref{eqn: h(0)= theta+ zg}$ and \eqref{g perp theta Hardy} together imply
\begin{equation}\label{h(0)}
h(0)=\beta \overline{\theta(0)}+\la zg, \theta \ra.
\end{equation}
Again, taking the inner product with $z^n\theta$ $(n\geq 1)$ on both sides of $\eqref{eqn: beta+g=theta h}$ and using \eqref{g perp theta Hardy}, we will have
$\la h, z^n\ra=0$ for all $n\geq 1$. Hence $h$ is a constant function on $\mathbb{\D}$: $h(z)=h(0)$.
Now substituting $z=\alpha\bar{\beta}$ in $\eqref{eqn: beta+g=theta h}$, we have $\beta= \theta(\alpha\bar{\beta})h(\alpha\bar{\beta})$. This implies,
$\theta(\alpha\bar{\beta})$ is nonzero and $h(\alpha\bar{\beta})=\frac{\beta}{\theta(\alpha\bar{\beta})}=h(0)$. Therefore,
$h(z)=\frac{\beta}{\theta(\alpha\bar{\beta})}$ and equation $\eqref{eqn: beta+g=theta h}$ reduces to

\begin{equation}\label{(z-ab)g= h(theta-ab)}
(z-\alpha\bar{\beta})g =\frac{\beta}{\theta(\alpha\bar{\beta})}\big(\theta(z)-\theta(\alpha\bar{\beta})\big).
\end{equation}
Hence $f=f_0+z^2g$, where $g$ is given by $\eqref{(z-ab)g= h(theta-ab)}$. Note that, $g$ depends only on $\theta$ and uniquely determined by the equation $\eqref{(z-ab)g= h(theta-ab)}$. Hence it follows that, $\dim \big(\clm\ominus z^2\theta H^2(\D)\big)=1$. Assume there exists another inner function $\theta_1$ such that $z^2\theta_1 H^2(\D)\subseteq\clm$ and $\dim \big(\clm\ominus z^2\theta_1 H^2(\D)\big)=1$. Then by construction of $\clm_1$ (see \eqref{M1}), $z^2\theta_1 H^2(\D)\subseteq\clm_1$ and hence by \eqref{M1 first}, $\theta$ divides $\theta_1$. Since $(\clm\ominus\clm_1)\big(=\clm\ominus z^2\theta H^2(\D)\big)$ is also of dimension one,  we must have $\theta_1= c\theta$ for some scalar $c$ with $|c|=1$. This completes the proof.
\end{proof}
Along the lines of proof of Theorem \ref{dim 1}, we would like to make the following remarks:
\begin{Remark}\label{rmk: g theta}
\begin{enumerate}
\item If $\clm$ is $S_{\alpha,\beta}$-invariant and $\dim\big(\clm\ominus z^2\theta H^2(\D)\big)=1$, then $\theta(\alpha\bar{\beta})\neq 0$.
\item An element $f \in \big(\clm\ominus z^2\theta H^2(\D)\big)$ can be written as $f=f_0+z^2g$, where $g\in (H^2(\D)\ominus
\theta H^2(\D))$ and satisfies $(z-\alpha\bar{\beta})g=\frac{\beta}{\theta(\alpha\bar{\beta})}(\theta(z)-\theta(\alpha\bar{\beta}))$.
\end{enumerate}
\end{Remark}

We now establish the cyclicity of $S_{\alpha,\beta}$-invariant subspaces. The following theorem also provides a characterization of the corresponding invariant subspaces.
\begin{Theorem}\label{Isp char}
Let $\clm\subseteq H^2_{\alpha,\beta}$ be a nonzero closed subspace. $S_{\alpha,\beta}\clm\subseteq\clm$ if and only if there exists an inner function $\theta$ such that either $\clm=z^2\theta H^2(\D)$ or $\clm$ can be written as
$$\clm=\C f \oplus z^2\theta H^2(\D), \text{ where}$$
$$f=f_0+ z^2g \text{ for some } g\in (H^2(\D)\ominus \theta H^2(\D)) \text{ satisfying}$$ $$(z-\alpha\bar{\beta})g=\frac{\beta}{\theta(\alpha\bar{\beta})}(\theta(z)-\theta(\alpha\bar{\beta})).$$ Moreover $\clm=[f]_{S_{\alpha,\beta}}$ i.e., $\clm$ is cyclic.
\end{Theorem}

\begin{proof}
Let $\clm\subseteq H^2_{\alpha,\beta}$ be a nonzero closed $S_{\alpha,\beta}$-invariant subspace. It follows by Theorem \ref{dim 1}, there exists an inner function $\theta$ such that, either $\clm= z^2\theta H^2(\D)$ or $\dim \big(\clm\ominus z^2\theta H^2(\D)\big)=1$. Also, an element $f\in\big(\clm\ominus z^2\theta H^2(\D)\big)$ can be written as $f= f_0+z^2g$, for some $g\in \big( H^2(\D)\ominus\theta H^2(\D)\big)$ such that
$$(z-\alpha\bar{\beta})g=\frac{\beta}{\theta(\alpha\bar{\beta})}[\theta(z)-\theta(\alpha\bar{\beta})].$$
Hence, we can write
$$\clm=\C f\oplus z^2\theta H^2(\D),$$ where $f$ satisfies the conditions as stated in the theorem.

Conversely, let $\clm\subseteq H^2(\D)$ be a nonzero closed subspace satisfying the given conditions. If $\clm=z^2\theta H^2(\D)$ for some inner function $\theta$, then $S_{\alpha,\beta}\clm=S_{\alpha,\beta}\big(z^2\theta H^2(\D)\big)=z^3\theta H^2(\D)\subseteq z^2\theta H^2(\D)=\clm$ and hence, $\clm$ is $S_{\alpha,\beta}$ invariant.

If $\clm=\C f\oplus z^2\theta H^2(\D)$ for some inner function $\theta$ and $f=f_0+ z^2g$ with $g\in \big(H^2(\D)\ominus \theta H^2(\D)\big)$ satisfying $(z-\alpha\bar{\beta})g=\frac{\beta}{\theta(\alpha\bar{\beta})}[\theta(z)-\theta(\alpha\bar{\beta})],$ we have
\begin{equation}\label{beta+ z g= theta}
\beta+ (z-\alpha\bar{\beta})g(z)=\frac{\beta}{\theta(\alpha\bar{\beta})}\theta(z)
\end{equation}
By a similar computation as in Theorem \ref{dim 1}, we have
\begin{equation}\label{theta g z}
S_{\alpha,\beta}f-\alpha\bar{\beta}f = S_{\alpha,\beta}(f_0+z^2g)-\alpha\bar{\beta} (f_0+z^2g)
= z^2(\beta+(z-\alpha\bar{\beta})g),
\end{equation}
which together with \eqref{beta+ z g= theta} imply $S_{\alpha,\beta}f=\alpha\bar{\beta}f+\frac{\beta}{\theta(\alpha\bar{\beta})}z^2\theta\in\clm$.
Also for any $h\in H^2(\D)$, $S_{\alpha,\beta}(z^2\theta h)= z^3\theta h\in z^2\theta H^2(\D)\subseteq\clm$. Hence $S_{\alpha,\beta}\clm\subseteq\clm.$

We now show $[f]_{S_{\alpha,\beta}}=\clm$. Note that,
\begin{equation}\label{cyclic [f]}
[f]_{S_{\alpha,\beta}} = \overline{span}\{f, S_{\alpha,\beta}f, S_{\alpha,\beta}^2f,\ldots\} =\overline{span}\{f, S_{\alpha,\beta}f-\alpha\bar{\beta}f, S_{\alpha,\beta}^2-\alpha\bar{\beta}S_{\alpha,\beta}f_1,\ldots\}=\C f+[S_{\alpha,\beta}f-\alpha\bar{\beta}f]_{S_{\alpha,\beta}}
\end{equation}
Again, it follows by \eqref{beta+ z g= theta}, \eqref{theta g z}, \eqref{eqn: S_ab matrix}, $$S_{\alpha,\beta}^n(S_{\alpha,\beta}f-\alpha\bar{\beta}f)=S^n(S_{\alpha,\beta}f-\alpha\bar{\beta}f)= S^n\Big(\frac{\beta}{\theta(\alpha\bar{\beta})}z^2\theta\Big), \text{ for all } n\geq 1.$$
Hence,
\begin{equation}\label{cyclic[Sf]}
[S_{\alpha,\beta}f-\alpha\bar{\beta}f]_{S_{\alpha,\beta}}=[S_{\alpha,\beta}f-\alpha\bar{\beta}f]_{S} = z^2\theta H^2(\D).
\end{equation}
Since $\clm=\C f\oplus z^2\theta H^2(\D)$, it follows by \eqref{cyclic [f]} and \eqref{cyclic[Sf]}, $[f]_{S_{\alpha,\beta}}=\C f\oplus z^2\theta H^2(\D)=\clm$.
\end{proof}
Note that, the equation \eqref{beta+ z g= theta} gives an implicit representation of the function $g$ associated to a generating element $f$ of the corresponding $S_{\alpha,\beta}$-invariant subspace. However if the inner function $\theta$ in \eqref{beta+ z g= theta} has a simple form, one can determine the function $g$ explicitly e.g., if $\theta$ is a monomial, the following corollary will show that $g$ is a polynomial.
\begin{Corollary}\label{cor1}
 In the setting of Theorem \ref{Isp char}, assume $\theta(z)=z^n$ for some $n\geq 1$. Then $\clm=[f]_{S_{\alpha,\beta}}$, where $f= f_0+ z^2g$, $g$ being the polynomial
$$g(z)=\frac{\beta}{(\alpha\bar{\beta})^n}[z^{n-1}+(\alpha\bar{\beta})z^{n-2}+(\alpha\bar{\beta})^2 z^{n-3}+\cdots+(\alpha\bar{\beta})^{n-2}z+(\alpha\bar{\beta})^{n-1}].$$
\end{Corollary}
\begin{proof}
In the proof of Theorem \ref{Isp char} above, substituting $\theta(z)=z^n$ $(n\geq 1)$ in \eqref{beta+ z g= theta}, we have
$$(z-\alpha\bar{\beta})g=\frac{\beta}{(\alpha\bar{\beta})^n}[z^n-(\alpha\bar{\beta})^n].$$
The proof is now immediate by unique factorization theorem for complex polynomials \cite{Artin}.
\end{proof}
In the setting of Theorem \ref{Isp char} and Corollary \ref{cor1}, We have the following observations:
\begin{enumerate}
\item $\clm\subseteq H^2_{\alpha,\beta}$ is of finite codimension if and only if $\theta$ is a finite Blaschke product.

\item If $\clm=[f]_{S_{\alpha,\beta}}$, $f$ need not be an inner function. For if $\theta(z)=z^n$ for some $n\geq1$, then by Corollary \ref{cor1},
$f=f_0+ \frac{\beta z^2}{(\alpha\bar{\beta})^n}[z^{n-1}+\cdots+(\alpha\bar{\beta})^{n-2}z+(\alpha\bar{\beta})^{n-1}]$, which is neither a constant nor a monomial with $|f(z)|=1$ on $|z|=1$.

\item If $\clm$, $\clm'$ are $S_{\alpha,\beta}$-invariant subspaces of $H^2_{\alpha,\beta}$ for some fixed pair $\{\alpha,\beta\}$, then
$S_{\alpha,\beta}|_{\clm}$ and
$S_{\alpha,\beta}|_{\clm'}$ need not be unitarily equivalent. Indeed, if $\clm= [z^2]_{S_{\alpha,\beta}}$, then $S_{\alpha,\beta}|_{\clm}$ is the unilateral shift of multiplicity one. If we take $\clm'= \C f\oplus z^2\theta H^2(\D)$ where $\theta(z)=z$, then by Corollary \ref{cor1}, $f$ can be taken as $f= f_0+\frac{\beta}{\alpha\bar{\beta}}z^2$. Then  $\{\frac{f}{\delta}, z^{n+3}, n\geq 0\}$ will form an orthonormal basis of $\clm'$, where $\delta=(1+\frac{1}{|\alpha|^2})^{1/2}$ and an easy calculation will show $\|S_{\alpha,\beta}(\frac{f}{\delta})\|=(|\alpha\beta|^2+\frac{1}{1+|\alpha|^2})^{1/2}>1$. Hence $S_{\alpha,\beta}|\clm'$ is not even an isometry.
\end{enumerate}
Note that, the final observation is in sharp contrast with a well-known consequence of the Beurling theorem: If $\clm_1$ and $\clm_2$ are nonzero closed $S$-invariant subspaces of $H^2(\D)$, then $S|_{\clm_1}$ and $S|_{\clm_2}$ are unitarily equivalent. In view of the observation $(3)$ above, this property fails to hold for $S_{\alpha,\beta}$-invariant subspaces.

Theorem \ref{dim 1} is not  quite an analogue of Beurling's Theorem:\textsf{ If $\cls$ is a nonzero closed $M_z$-invariant subspace of $H^2(\D)$, then $\dim(\cls\ominus z\cls)=1$.} In the next section, we prove a similar result for $S_{\alpha,\beta}$-invariant subspaces. Based on that, we provide another characterization of $S_{\alpha,\beta}$-invariant subspaces which will be useful to discuss their wandering subspace property.

\newsection{Beurling type Theorem and Invariant Subspaces}\label{sec: Characterization}

\begin{Theorem}\label{codim 1}
If $\clm\subseteq H^2_{\alpha,\beta}$ is a nonzero, closed, $S_{\alpha,\beta}$-invariant subspace, then
\[\dim(\clm\ominus S_{\alpha,\beta}\clm)=1.\]
\end{Theorem}
\begin{proof}
If $\clm=z^2\theta H^2(\D)$ for some inner function $\theta$, then $$S_{\alpha,\beta}\clm=S_{\alpha,\beta}(z^2\theta H^2(\D))=S(z^2\theta H^2(\D))=z^3\theta H^2(\D)$$ and hence
$\dim(\clm\ominus S_{\alpha,\beta}\clm)=1$. So let us assume there exists $h\in\clm$ such that $h(0)\neq 0$. Since $\clm$ is closed and $S_{\alpha,\beta}$ is left invertible (see Theorem \ref{shift}), $S_{\alpha,\beta}\clm$ is also closed in $H^2_{\alpha,\beta}$. Again $S_{\alpha,\beta}\clm\subseteq\clm$ implies $S_{\alpha,\beta}(S_{\alpha,\beta}\clm)\subseteq S_{\alpha,\beta}\clm$ i.e., $S_{\alpha,\beta}\clm$ is a closed $S_{\alpha,\beta}$-invariant subspace of $H^2_{\alpha,\beta}$. We now split the proof into the following steps.

\NI\textsf{Step 1:} In this step, we will provide a representation of $S_{\alpha,\beta}\clm$.

\NI Since $h\in\clm$ and $h(0)\neq 0$, w.l.g. we can write $S_{\alpha,\beta}h= S_{\alpha,\beta}(f_0+z^2h_1)$, for some $h_1\in H^2(\D)$. Then (by \eqref{eqn: S_ab matrix}), $S_{\alpha,\beta}h= \alpha\bar{\beta} f_0+\beta z^2+z^3h_1\in S_{\alpha,\beta}\clm$ with $S_{\alpha,\beta}h(0)=\alpha^2\bar{\beta}\neq 0$. Hence by Theorem \ref{Isp char}, there exists an inner function $\theta$ such that
\begin{eqnarray}
S_{\alpha,\beta}\clm &=& \C f\oplus z^2\theta H^2(\D), \text{ where} \label{M}\\
f &=& f_0+z^2 g,\ \ g\in \big(H^2(\D)\ominus \theta H^2(\D)\big)\text{ and } \label{f}\\
(z-\alpha\bar{\beta})g &=&\frac{\beta}{\theta(\alpha\bar{\beta})}\big(\theta(z)-\theta(\alpha\bar{\beta})\big) \label{g new}.
\end{eqnarray}
Note that, along the lines of proof of the Theorem \ref{dim 1}, the inner function $\theta$ corresponding to an $\clm$ is precisely associated with the subspace $\overline{span}\{f\in \clm: f(z)= z^2g(z)\text{ for some } g\in H^2(\D)\}= z^2\theta H^2(\D)$ ( see \eqref{M1}, \eqref{M1 first}).
In the setting of $S_{\alpha,\beta}\clm$, this reduces to the following one (which we call by $\clm_1$ again):
\begin{equation}\label{S(M1) modified}
\clm_1= \overline{span}\{S_{\alpha,\beta}f\in S_{\alpha,\beta}\clm: S_{\alpha,\beta}f(z)= z^2g(z)\text{ for some } g\in H^2(\D)\}= z^2\theta H^2(\D).
\end{equation}
Note, corresponding to $h$ above, $S_{\alpha,\beta}^2h-\alpha\bar{\beta} S_{\alpha,\beta}h=z^3\big(\beta+(z-\alpha\bar{\beta})h_1\big)\in\clm_1$. We show $z$ divides $\theta$. On the contrary, suppose $z\not|\theta$. Then there exists $f\in\clm_1$ such that $f(z)=z^2 g(z)$ for some $g\in H^2(\D)$ with $g(0)\neq 0$. Also, there exist $c_0\in \C$ and $h_2\in H^2(\D)$ such that
\begin{equation}\label{f B-type}
z^2 g=f= S_{\alpha,\beta}(c_0f_0+z^2h_2)=c_0(\alpha\bar{\beta}f_0+\beta z^2)+z^3h_2
\end{equation}
Since $g(0)\neq 0$, \eqref{f B-type} implies $c_0\neq 0$. Again on substitution of $z=0$, \eqref{f B-type} yields, $0=c_0(\alpha^2\beta)$ i.e., $c_0=0$, which is a contradiction. Hence there exists an inner function $\theta_1$ such that
\begin{eqnarray}
\theta(z) &=& z\theta_1(z) \label{theta=ztheta1}\\
\clm_1 &=& z^3\theta_1 H^2(\D)\ \ (\text{ by \ref{S(M1) modified}})\label{M1=z3 theta1}),
\end{eqnarray}
and hence, \eqref{S(M1) modified} can be rewritten as
\begin{equation}\label{S(M1) modified to 3}
\clm_1= \overline{span}\{S_{\alpha,\beta}f\in S_{\alpha,\beta}\clm: S_{\alpha,\beta}f(z)= z^3g(z)\text{ for some } g\in H^2(\D)\}= z^3\theta_1 H^2(\D).
\end{equation}
Note that, by Remark \ref{rmk: g theta}, $\theta(\alpha\bar{\beta})\neq 0$. Hence it follows by \eqref{theta=ztheta1}, $\theta_1(\alpha\bar{\beta})\neq 0$. Since $\dim(S_{\alpha,\beta}\clm\ominus \clm_1)=1$ (by \eqref{M} and \eqref{S(M1) modified}), up to a scalar multiple, an element $f_1\in~(S_{\alpha,\beta}\clm\ominus \clm_1)$ can be written as $f_1= S_{\alpha,\beta}(f_0+z^2 h_3)$ for some $h_3\in H^2(\D)$ and hence
\begin{equation}\label{f prime for S1M}
f_1= \alpha\bar{\beta}f_0+ z^2(\beta+z h_3).
\end{equation}
Also, note that
\begin{equation}\label{f prime /ab}
S_{\alpha,\beta}\clm\ominus\clm_1= \C f_1=\C \big(\frac{f_1}{\alpha\bar{\beta}}\big).
\end{equation}
By \eqref{f} and \eqref{f prime for S1M}, $f$ can be taken as
\begin{equation}\label{f modified}
f=\frac{1}{\alpha\bar{\beta}}f_1= f_0+\frac{z^2}{\alpha\bar{\beta}}(\beta+ zh_3).
\end{equation}
Again by \eqref{f} and \eqref{f modified}, the corresponding function $g$ will be
\begin{equation}\label{g modified}
g(z)=\frac{1}{\alpha\bar{\beta}}(\beta+ zh_3).
\end{equation}
We now show
\begin{equation}\label{ip (h3, theta)=0}
 h_3\in \big(H^2(\D)\ominus \theta_1 H^2(\D)\big).
\end{equation}
Indeed, $g\perp\theta H^2(\D)$ implies $\la g, z^n\theta\ra=0$ for all $n\geq 0$. This together with \eqref{g modified} and \eqref{theta=ztheta1} implies $\la\beta+ zh_3, z^{n+1}\theta_1\ra=0$ and hence $\la h_3,z^n\theta_1\ra=0$ for all $n\geq 0$. Again by \eqref{g new} and \eqref{g modified},
\begin{equation}\label{h3 calculation}
\frac{1}{\alpha\bar{\beta}}(z-\alpha\bar{\beta})(\beta+z h_3)=\frac{\beta}{\theta(\alpha\bar{\beta})}\Big(\theta(z)-\theta(\alpha\bar{\beta})\Big)=\frac{\beta}{\alpha\bar{\beta}\theta_1(\alpha\bar{\beta})}\Big(z\theta_1(z)-\alpha\bar{\beta} \theta_1(\alpha\bar{\beta})\Big),
\end{equation}
where the last equality follows by \eqref{theta=ztheta1}. A simple computation by equating the left most term with that of the right most term of \eqref{h3 calculation} yields,
\begin{equation}\label{z involve h3}
z\big(\beta+(z-\alpha\bar{\beta})h_3\big)=\frac{\beta}{\theta_1(\alpha\bar{\beta})}z\theta_1.
\end{equation}
Since $\beta+(z-\alpha\bar{\beta})h_3\neq 0$, the equation \eqref{z involve h3} reduces to
\begin{equation}\label{final h3}
(z-\alpha\bar{\beta})h_3=\frac{\beta}{\theta_1(\alpha\bar{\beta})}\big(\theta_1(z)-\theta_1(\alpha\bar{\beta})\big).
\end{equation}
Now, the equations \eqref{M1=z3 theta1}, \eqref{f prime /ab}, \eqref{f modified}, and \eqref{ip (h3, theta)=0} together imply
$$S_{\alpha,\beta}\clm=\C f_1\oplus z^3\theta_1 H^2(\D),\quad \text{where}$$
$f_1=\alpha\bar{\beta}f_0+\beta z^2+ z^3h_3$ with $h_3\in H^2(\D)\ominus\theta_1H^2(\D)$ such that \eqref{final h3} holds.
Denoting (by abuse of notation, as there will be no harm) $\theta_1, h_3$ by $\theta$ and $g$ respectively, it summarizes to
\begin{eqnarray}
S_{\alpha,\beta}\clm &=& \C f_1\oplus z^3\theta H^2(\D),\quad \text{where} \label{S1M last}\\
f_1 &=& \alpha\bar{\beta}f_0+\beta z^2+ z^3g, \quad g\in H^2(\D) \text{ such that} \label{f1 last}\\
\la g, z^n\theta\ra &=& 0,  \text{ for all } n\geq 0 \text{ and} \label{g last}\\
(z-\alpha\bar{\beta})g &=& \frac{\beta}{\theta(\alpha\bar{\beta})}\big(\theta(z)-\theta(\alpha\bar{\beta})\big). \label{g equation}
\end{eqnarray}

\NI\textsf{Step 2:} We now show that, for any nonzero subspace $\clm$ of $H^2_{\alpha,\beta}$, $S_{\alpha,\beta}\clm\neq\clm$.

\NI Let $\clm= H^2_{\alpha,\beta}$ and suppose $S_{\alpha,\beta}\clm=\clm$. Then $S_{\alpha,\beta}$ is invertible as $S_{\alpha,\beta}$ is left invertible (by Theorem \ref{shift}). This leads to a contradiction, as by Lemma \ref{kernel S_ab adjoint}, $\ker S_{\alpha,\beta}^*\neq \{0\}$. Hence $S_{\alpha,\beta} H^2_{\alpha,\beta}\neq H^2_{\alpha,\beta}$.
Now let $\clm$ be a nonzero, proper, closed subspace of $H^2_{\alpha,\beta}$ such that $S_{\alpha,\beta}\clm=\clm$. This yields, $S_{\alpha,\beta}^n\clm=\clm$ for all $n\geq 1$ and hence
$$\clm\subseteq\bigcap_{n\geq 1}S_{\alpha,\beta}^n H^2_{\alpha,\beta}=\{0\},$$

where the last equality follows by the analyticity of $S_{\alpha,\beta}$ (Theorem \ref{shift}). This contradicts our assumption $\clm\neq \{0\}$. Hence for any nonzero closed subspace $\clm\subseteq H^2_{\alpha,\beta}$, $S_{\alpha,\beta}\clm\neq \clm$.

\NI\textsf{Step 3:} In this step, we provide a representation of a generating element of $\clm\ominus S_{\alpha,\beta}\clm$. This generating function turns out to be unique upto a scalar multiple.

\NI By \textsf{Step 2}, $\exists$ $f_2(\neq 0)\in\clm\ominus S_{\alpha,\beta}\clm$. Up to a nonzero scalar multiple, we can write $f_2$ as
\begin{equation}\label{f2 defined}
f_2=f_0+z^2g_1, \text{ for some } g_1\in H^2(\D).
\end{equation}
Proceeding exactly in the same way as before, we will have
\begin{eqnarray}
S_{\alpha,\beta}f_2 &=& (\alpha\bar{\beta}f_0+\beta z^2)+z^3g_1\ \ \label{f 2}\\
S_{\alpha,\beta}^2f_2-\alpha\bar{\beta}S_{\alpha,\beta}f_2 &=& z^3[\beta+(z-\alpha\bar{\beta})g_1](\neq 0)\in \clm_1=z^3\theta H^2(\D)\quad ( \text{ by } \eqref{S(M1) modified to 3})\label{f2=h'' first}
\end{eqnarray}
Now, $S_{\alpha,\beta}f_2\in S_{\alpha,\beta}\clm$. Hence by \eqref{S1M last}, there exists $h_4\in H^2(\D)$ and $c_0\in\C$ such $S_{\alpha,\beta}f_2= c_0f_1+z^3\theta h_4$.
This together with \eqref{f1 last} implies,
\begin{equation}\label{s1f2}
\begin{split}
S_{\alpha,\beta}f_2= c_0f_1+z^3\theta h_4
&=c_0(\alpha\bar{\beta}f_0+\beta z^2+z^3 g)+ z^3\theta h_4
\end{split}
\end{equation}
Hence it follows by \eqref{f 2} and \eqref{s1f2}, $c_0=1$ and
\begin{equation}\label{z3g_1=z3(g+theta h')}
z^3 g_1 = z^3(g+\theta h_4)
\end{equation}
By \eqref{f2=h'' first}, there exists $h_5\in H^2(\D)$ such that,
\[
S_{\alpha,\beta}^2f_2-\alpha\bar{\beta}S_{\alpha,\beta}f_2=z^3[\beta+(z-\alpha\bar{\beta})g_1]= z^3\theta h_5 \text{ and hence}
\]
\begin{equation}\label{beta+(z-ab)g_1=theta h''}
\beta+(z-\alpha\bar{\beta})g_1 =\theta h_5\quad (\text{ as } \beta+(z-\alpha\bar{\beta})g_1\neq 0).
\end{equation}
Since $f_2\in\clm\ominus S_{\alpha,\beta}\clm$, by \eqref{S1M last}
\begin{eqnarray}
\la f_2, f_1\ra&=& 0 \label{ip (f1,f2)=0}, \text{ and}\\
\la f_2, z^{n+3}\theta \ra &=& 0\ \ \text{ for all } n\geq 0. \label{ip (f2, zn theta)=0}
\end{eqnarray}
By \eqref{f2 defined} and \eqref{ip (f2, zn theta)=0}, $\la f_0+z^2 g_1, z^{n+3}\theta\ra=0$ and hence,
\begin{equation}\label{g1 perp on z theta Hardy space}
\la g_1,z^{n+1}\theta\ra=0 \text{ for all } n\geq 0.
\end{equation}
Again \eqref{f1 last}, \eqref{f2 defined}, and \eqref{ip (f1,f2)=0} together imply, $$\la f_0+z^2g_1, \alpha\bar{\beta}f_0+z^3g+\beta z^2\ra =0,$$ and a further simplification reveals $\bar{\alpha}\beta +\la g_1,zg\ra+ \bar{\beta}\la g_1, 1\ra =0$. This together with \eqref{g equation} (substituting $zg$) implies
\begin{equation}\label{f2 f1 zero i.p}
\bar{\alpha}\beta +\bar{\alpha}\beta\la g_1,g\ra+\frac{\bar{\beta}}{\overline{\theta(\alpha\bar{\beta})}}\la g_1,\theta\ra =0
\end{equation}
Note that, $g_1\neq 0$. Otherwise, \eqref{f2 f1 zero i.p} would imply $\bar{\alpha}\beta=0$. Hence by \eqref{z3g_1=z3(g+theta h')},
\begin{equation}\label{g1=g+theta h'}
g_1= g+\theta h_4.
\end{equation}
Note that, $g, g_1$ are both orthogonal to $z^n\theta H^2(\D)$ for $n\geq 1$ ( by \eqref{g last}, \eqref{g1 perp on z theta Hardy space}). Then, the inner product with $z^n\theta$ ( $n\geq 1$) on both sides of \eqref{g1=g+theta h'} yields, $\la h_4, z^n\ra=0$.
Hence \textsf{ $h_4$ is a constant.} Again, as $g\perp \theta H^2(\D)$ \eqref{g last}, we have by \eqref{g1=g+theta h'},
\begin{equation}\label{h'(0)}
\la g_1,\theta\ra =\la h_4,1\ra=h_4(0).
\end{equation}
We show  $\la g_1,\theta\ra\neq 0$. On the contrary, let $\la g_1,\theta\ra=0$. Then \eqref{h'(0)} implies, $h_4(z)=h_4(0)=0$ and further by \eqref{g1=g+theta h'}, we have $g_1=g$. Then by \eqref{f2 f1 zero i.p} and \eqref{g last}, $(1+\|g\|^2)=0$, which is impossible.
Hence by \eqref{h'(0)}, $h_4(z)=h_4(0)=\la g_1,\theta\ra\neq 0,$ and \eqref{g1=g+theta h'} implies
\begin{equation}\label{g1=g+h'theta}
g_1(z)=g(z)+h_4(0)\theta(z).
\end{equation}
Now substituting $g_1$ ( from\eqref{g1=g+h'theta}) in the equation \eqref{f2 f1 zero i.p} and using \eqref{g last}, one will have
\begin{equation}\label{h'(0)=}
h_4(0) = -\frac{\bar{\alpha}\beta\overline{\theta(\alpha\bar{\beta})} (1+\|g\|^2)}{\bar{\beta}}.
\end{equation}
Hence, by \eqref{g1=g+h'theta} and \eqref{h'(0)=}, the function $g_1$ can be written as
\begin{equation}\label{g1=g-r theta}
g_1(z)=g(z)-\frac{\bar{\alpha}\beta\overline{\theta(\alpha\bar{\beta})} (1+\|g\|^2)}{\bar{\beta}}\theta(z)
\end{equation}
Hence by \eqref{f2 defined} and\eqref{g1=g-r theta}, one can write $f_2$ as
\begin{equation}\label{f2}
f_2 = f_0+z^2g-\frac{\bar{\alpha}\beta}{\bar{\beta}}\overline{\theta(\alpha\bar{\beta})}(1+\|g\|^2)z^2\theta.
\end{equation}
Since (as we noted earlier), $g$ is uniquely determined by $\theta$, the same is true for $g_1$ (by \eqref{g1=g-r theta}. Hence by \eqref{f2}, $f_2$ is uniquely determined by $\theta$. This is equivalent to saying that, any function in $\big(\clm\ominus S_{\alpha,\beta}\clm\big)$ can differ from $f_2$ only by a scalar multiple. Hence, $\dim(\clm\ominus S_{\alpha,\beta}\clm)=1$.
\end{proof}
\begin{Remark}
 We have used Theorem \ref{dim 1} to represent $S_{\alpha,\beta}\clm$ in Theorem \ref{codim 1}.
\end{Remark}
In light of Theorem \ref{codim 1}, we have the following characterization result.
\begin{Theorem}\label{S1 inv 1}
Let $\clm\subseteq H^2_{\alpha,\beta}$ be a nonzero closed subspace. Then $S_{\alpha,\beta}\clm\subseteq\clm$ if and only if there exists an inner function $\theta$ and a function $g\in \big(H^2(\D)\ominus\theta H^2(\D)\big)$ such that
\begin{equation*}
\begin{split}
\clm &=\mathbb{C} f_2\oplus\mathbb{C} f_1\oplus z^3\theta H^2(\D), \text{ where}\\
f_1 &=\alpha\bar{\beta}f_0+\beta z^2+z^3g, \text{ with  } \\
(z-\alpha\bar{\beta})g &=\frac{\beta}{\theta(\alpha\bar{\beta})}[\theta(z)-\theta(\alpha\bar{\beta})], \text{ and} \\
f_2 &= f_0+z^2g-\frac{\bar{\alpha}\beta}{\bar{\beta}}\overline{\theta(\alpha\bar{\beta})}(1+\|g\|^2)z^2\theta.
\end{split}
\end{equation*}
\end{Theorem}
\begin{proof}
Let $\clm\subseteq H^2_{\alpha,\beta}$ be a nonzero closed $S_{\alpha,\beta}$- invariant subspace. Then following the lines of proof of Theorem \eqref{codim 1}, there exist $\theta$ and $g$ such that $\clm$ can be written as in the statement.

Conversely, let $\clm\subseteq H^2_{\alpha,\beta}$ be a nonzero closed subspace satisfying the given conditions of the theorem. We show that,
$S_{\alpha,\beta}\clm\subseteq\clm$.
\\
Note that (by \eqref{eqn: S_ab matrix}), for any $h\in H^2(\D)$, $S_{\alpha,\beta}(z^3\theta h)=z^4\theta h\subseteq z^3\theta H^2(\D)$. Next,
\begin{equation}\label{s1f2=f1-theta}
\begin{split}
S_{\alpha,\beta}f_2 &= S_{\alpha,\beta}\Big(f_0+z^2g-\frac{\bar{\alpha}\beta}{\bar{\beta}}\overline{\theta(\alpha\bar{\beta})}(1+\|g\|^2)z^2\theta\Big)= (\alpha\bar{\beta}f_0+\beta z^2)+z^3g-\frac{\bar{\alpha}\beta}{\bar{\beta}}\overline{\theta(\alpha\bar{\beta})}(1+\|g\|^2)z^3\theta\\ &=f_1-\frac{\bar{\alpha}\beta}{\bar{\beta}}\overline{\theta(\alpha\bar{\beta})}(1+\|g\|^2)z^3\theta \in\Big(\C f_1\oplus z^3\theta H^2(\D)\Big)\subseteq\clm,
\end{split}
\end{equation}
where the third equality follows by the expression of $f_1$ given in the statement. Now by \eqref{eqn: S_ab matrix},
\begin{eqnarray}
S_{\alpha,\beta}f_1&=& S_{\alpha,\beta}(\alpha\bar{\beta}f_0+z^3g+\beta z^2)=\alpha\bar{\beta}(\alpha\bar{\beta}f_0+\beta z^2)+z^4g+\beta z^3, \text{ and} \\
S_{\alpha,\beta}f_1-\alpha\bar{\beta}f_1 &=& S_{\alpha,\beta}f_1-\alpha\bar{\beta}(\alpha\bar{\beta}f_0+\beta z^2+z^3g)=z^3\big(\beta+(z-\alpha\bar{\beta})g\big) \label{S1 f1-ab f1}.
\end{eqnarray}
Again, the given condition $(z-\alpha\bar{\beta})g =\frac{\beta}{\theta(\alpha\bar{\beta})}[\theta(z)-\theta(\alpha\bar{\beta})]$ together with \eqref{S1 f1-ab f1} implies,
\begin{equation*}
S_{\alpha,\beta}f_1-\alpha\bar{\beta}f_1=\frac{\beta}{\theta(\alpha\bar{\beta})}z^3\theta.
\end{equation*}
Hence, $S_{\alpha,\beta}f_1 =\alpha\bar{\beta}f_1+\frac{\beta}{\theta(\alpha\bar{\beta})}z^3\theta\in\C f_1\oplus z^3\theta H^2(\D)\subseteq\clm $ and we conclude $S_{\alpha,\beta}\clm\subseteq\clm$.
\end{proof}
\begin{Remark}
\begin{enumerate}
\item The observations following Lemma \ref{cor1} on $S_{\alpha,\beta}$-invariant subspace, will remain same in the setting of Theorem \ref{S1 inv 1}.
\item The equations \eqref{s1f2=f1-theta} and \eqref{S1 f1-ab f1} together imply, the function $f_2$ in Theorem \ref{S1 inv 1} satisfies the property
$f_2\perp S_{\alpha,\beta}^nf_2$ for all $n\in\N$.
\end{enumerate}
\end{Remark}

\newsection{Wandering Subspace property}\label{sec: Wandering Subspace property}

In this section, we discuss the wandering subspace property (w.s.p.) of a nonzero closed $S_{\alpha,\beta}$-invariant subspace $\clm\subseteq H^2_{\alpha,\beta}$. We have seen in the earlier section (Theorem \ref{codim 1}), $\dim(\clm\ominus S_{\alpha,\beta}\clm)=1$. In Theorem \ref{Isp char}, we establish the cyclicity of $S_{\alpha,\beta}$-invariant subspaces. At this point, one may wish to investigate if $[\clm\ominus S_{\alpha,\beta}\clm]_{S_{\alpha,\beta}}= \clm$  holds just like other well-known classical shift operators. However, this is not so. In this section, we provide counter examples and discuss to which extent w.s.p. of an $S_{\alpha,\beta}$-invariant subspace can hold.

To begin with, first we consider the full space $H^2_{\alpha, \beta}$ and verify if $ H^2_{\alpha,\beta}=[H^2_{\alpha,\beta}\ominus S_{\alpha,\beta} H^2_{\alpha,\beta}]_{S_{\alpha,\beta}}=[ker S_{\alpha,\beta}^*]_{S_{\alpha,\beta}}$. We start with the following lemma.
Let us denote:
\begin{equation}\label{t defined}
p=\frac{\bar{\beta}}{\bar{\alpha}}(1+|\alpha|^2).
\end{equation}
Clearly $p\neq 0$  as $\beta$ is nonzero. Then,
\begin{Lemma}\label{t<1}
$|p|\geq1$ if and only if $|\alpha|^2\leq \frac{1}{u+1}$, where $u$ is the unique real root of the cubic $z^3+3z^2+2z-1$.
\end{Lemma}
\begin{proof}
Substituting $|\beta|^2=1-|\alpha|^2$ in \eqref{t defined}, a simple computation shows, $|p|\geq 1$ is equivalent to $|\alpha|^6+|\alpha|^4-1\leq 0$. Let us denote $x:= |\alpha|^2$. Then we have,
\begin{equation}\label{alpha and x}
|\alpha|^6+|\alpha|^4-1\leq 0 \iff x^3+(x+1)(x-1)\leq 0.
\end{equation}
Since $x(=|\alpha|^2)<1$, we have $x=\frac{1}{1+y}$ for some $y>0$. Now, on substituting $x=\frac{1}{1+y}$, the relation $ x^3+(x+1)(x-1)\leq 0$ can be seen to be equivalent to $ y(1+y)(2+y)-1\geq 0$. Again, the expression $y(1+y)(2+y)-1$ is equivalent to $y^3+3y^2+2y-1$. Let us consider the function $f: \mathbb{R} \longrightarrow \mathbb{R}$ defined by $f(y)= y(1+y)(2+y)-1$. Note that, $f(0)=-1$, $f(1)=5$ and $f'(y)= 3 y^2+6y+2>0$ for $y\in [0,\infty)$. Hence the function $f$ is strictly increasing with exactly one root in $(0, \infty)$. Let us call this root $u$, so that $f(y)\geq 0$ if and only if $y\geq u$. It is easy to see that the other roots of $f$ are imaginary. Since $y\geq u$ is equivalent to $x\leq\frac{1}{1+u}$, the proof is complete by \eqref{alpha and x}.
\end{proof}
We noted earlier (Lemma \ref{kernel S_ab adjoint}), $\ker S_{\alpha,\beta}^*=\C(f_0-\frac{\bar{\alpha}\beta}{\bar{\beta}}z^2)$. Let us denote $\tilde{f}:=f_0-\frac{\bar{\alpha}\beta}{\bar{\beta}}z^2$. With this notion, we now prove the following theorem on w.s.p. of the full space $H^2_{\alpha,\beta}$:
\begin{Theorem}\label{alpha less zero}
For $\alpha,\beta\neq 0$ and $|\alpha|^2+|\beta|^2=1$, $H^2_{\alpha,\beta}=[H^2_{\alpha,\beta}\ominus S_{\alpha,\beta} H^2_{\alpha,\beta}]_{S_{\alpha,\beta}}$ if and only if $|\alpha|^2\leq \frac{1}{u+1}$, $u$ being the unique real root of the cubic polynomial $z^3+3z^2+2z-1$.
\end{Theorem}
\begin{proof}
We only need to verify $[\tilde{f}]_{S_{\alpha,\beta}}= H^2_{\alpha,\beta}$. Note that,
\begin{equation*}
S_{\alpha,\beta}\tilde{f}-\alpha\bar{\beta}\tilde{f} = S_{\alpha,\beta}(f_0-\frac{\bar{\alpha}\beta}{\bar{\beta}}z^2)-\alpha\bar{\beta}(f_0-\frac{\bar{\alpha}\beta}{\bar{\beta}}z^2)= (\alpha\bar{\beta}f_0+\beta z^2-\frac{\bar{\alpha}\beta}{\bar{\beta}}z^3)-\alpha\bar{\beta}(f_0-\frac{\bar{\alpha}\beta}{\bar{\beta}}z^2)\ \ (\text{by} \eqref{eqn: S_ab matrix}),
\end{equation*}
and hence
\begin{equation}\label{s1f1=t}
S_{\alpha,\beta}\tilde{f}-\alpha\bar{\beta}\tilde{f}=\frac{\bar{\alpha}\beta}{\bar{\beta}} z^2\Big(\frac{\bar{\beta}}{\bar{\alpha}}(1+|\alpha|^2)-z\Big)=\frac{\bar{\alpha}\beta}{\bar{\beta}} z^2(p-z),
\end{equation}
where the last equality holds by \eqref{t defined}. Since
\begin{equation*}
\begin{split}
[\tilde{f}]_{S_{\alpha,\beta}} =\overline{span}\{\tilde{f}, S_{\alpha,\beta}\tilde{f}, S_{\alpha,\beta}^2\tilde{f},\ldots\} &=\overline{span}\{\tilde{f}, S_{\alpha,\beta}\tilde{f}-\alpha\bar{\beta}\tilde{f}, S_{\alpha,\beta}^2\tilde{f}-\alpha\bar{\beta}S_{\alpha,\beta}\tilde{f},\ldots\}
\\
&=\C \tilde{f}+[S_{\alpha,\beta}\tilde{f}-\alpha\bar{\beta}\tilde{f}]_{S_{\alpha,\beta}},
\end{split}
\end{equation*}
it follows by \eqref{s1f1=t}, \eqref{t defined}, and \eqref{eqn: S_ab matrix}
\begin{equation}\label{[f1]=z theta}
[\tilde{f}]_{S_{\alpha,\beta}}=\C \tilde{f} +[z^2(p-z)]_{S_{\alpha,\beta}}=\C \tilde{f} +z^2[(p-z)]_S.
\end{equation}
Let $H^2_{\alpha, \beta}=[\tilde{f}]_{S_{\alpha,\beta}}$. Since $\tilde{f}(0)=f_0(0)=\alpha\neq 0$, by \ref{[f1]=z theta}, we must have $z^2\in z^2[p-z]_{S}$. If $|p|<1$, the function $f(z)=\sum_{n=0}^{\infty}\bar{p}^n z^n\in H^2(\D)$ and it is easy to check that $f\perp z^n(p-z)$ for all $n\geq 0$. Hence,
$f\perp [p-z]_S$ and there exists a non-constant inner function $\theta$ such that $z^2[p-z]_S=z^2\theta H^2(\D)$. Now $z^2\in z^2[p-z]_{S}\big(=z^2\theta H^2(\D)\big)$ implies $\theta$ must be a constant (of unit modulus), which is a contradiction. Hence $|p|\geq1$ and the proof of this part follows by the forward implication of Lemma \eqref{t<1}.

Conversely, let $|\alpha|^2\leq\frac{1}{1+u}$, as in the statement. Then by Lemma \eqref{t<1}, $|p|\geq1$ and hence $[p-z]_S=H^2(\D)$. This together with \eqref{[f1]=z theta} implies, $[\tilde{f}]_{S_{\alpha,\beta}}=\C \tilde{f} + z^2 H^2(\D)$. Since $\tilde{f}= (f_0-\frac{\bar{\alpha}\beta}{\bar{\beta}}z^2)$, it follows that, $f_0\in \C \tilde{f}+ z^2 H^2(\D)$ and hence $[\tilde{f}]_{S_{\alpha,\beta}}=H^2_{\alpha, \beta}$.
\end{proof}
Let us now turn to the nonzero closed and proper subspaces $\clm\subseteq H^2_{\alpha,\beta}$, invariant under $S_{\alpha,\beta}$. We have noticed in Theorem \ref{codim 1}, $(\clm \ominus S_{\alpha,\beta}\clm)=\C f_2$ where $f_2$ is a nonzero element of $\clm$ given by \eqref{f2}. Based on this, Theorem \ref{S1 inv 1} provides a complete characterization of $S_{\alpha,\beta}$-invariant subspaces. Our goal is now to verify whether $f_2$ generates $\clm$. Let us set
\begin{equation}\label{r}
r=\frac{\alpha\bar{\beta}(1+|\alpha|^2|\theta(\alpha\bar{\beta})|^2(1+\|g\|^2)}{|\alpha|^2|\theta(\alpha\bar{\beta})|^2(1+\|g\|^2)},
\end{equation}
where $\theta$ and $g$ are functions associated to $\clm$ (Theorem \eqref{S1 inv 1}). Recall that, by Theorem \ref{S1 inv 1}, there exists an inner function $\theta$ such that, $\clm$ can be represented as
\begin{eqnarray}
\clm &=& \C f_2\oplus\C f_1\oplus z^3\theta H^2(\D)\text{ where} \label{isp1} \\
f_1 &=& \alpha\bar{\beta}f_0+z^3g+\beta z^2, \text{ with }  g\in \big(H^2(\D)\ominus\theta H^2(\D)\big) \text{ satisfying}\label{isp2} \\
(z-\alpha\bar{\beta})g &=& \frac{\beta}{\theta(\alpha\bar{\beta})}[\theta(z)-\theta(\alpha\bar{\beta})], \text{ and}\label{isp3} \\
f_2&=& f_0+z^2g-\frac{\bar{\alpha}\beta}{\bar{\beta}}\overline{\theta(\alpha\bar{\beta})}(1+\|g\|^2)z^2\theta,\text{ with } f_2\perp S_{\alpha,\beta}\clm.\label{isp4}
\end{eqnarray}

Again, Recall by \eqref{f2=h'' first}, \eqref{g1=g+h'theta}, \eqref{h'(0)=}, and \eqref{beta+(z-ab)g_1=theta h''},
\begin{eqnarray}
S_{\alpha,\beta}^2f_2-\alpha\bar{\beta}S_{\alpha,\beta}f_2 &=& z^3\big(\beta+(z-\alpha\bar{\beta})g_1\big),\ \ \text{ where } \label{s12 f2 same} \\
g_1 &=& g(z)+h_4(0) \theta(z), \label{g1 recall}\\
h_4(z)=h_4(0)&=&-\frac{\bar{\alpha}\beta\overline{\theta(\alpha\bar{\beta})} (1+\|g\|^2)}{\bar{\beta}}\ \ \text{ and}\label{h4 recall}\\
\beta+(z-\alpha\bar{\beta})g_1 &=& \theta h_5,\ \ \text{ for some } h_5\in H^2(\D).\label{beta+(z-ab)g_1=theta h'' same}
\end{eqnarray}
At this point, we have the following lemma:
\begin{Lemma}\label{polynomial rep}
$h_5(z)$ is a linear polynomial and is given by
$$h_5(z)=\frac{\bar{\alpha}\beta\overline{\theta(\alpha\bar{\beta})}(1+\|g\|^2)}{\bar{\beta}} (r-z),$$
where $r$ is given by \eqref{r}.
\end{Lemma}
\begin{proof}
Taking inner product with $\theta$ on both sides of \eqref{beta+(z-ab)g_1=theta h'' same}, we have
\begin{equation}\label{wsp1}
\beta\overline{\theta(0)}+\la zg_1,\theta\ra-\alpha\bar{\beta}\la g_1,\theta\ra = h_5(0).
\end{equation}
Substituting $g_1$ ( from \eqref{g1 recall}) and using $g\perp \theta H^2(\D)$ (\eqref{g last}), the equation \eqref{wsp1} further reduces to
\begin{equation}\label{wsp2}
\beta\overline{\theta(0)}+\la zg,\theta\ra-\alpha\bar{\beta} h_4(0)= h_5(0).
\end{equation}
Again, substituting $zg$ (from \eqref{isp3}) in the above equation \eqref{wsp2} and using \eqref{g last}, a simple computation shows
\begin{equation}\label{wsp3}
 h_5(0)= \frac{\beta}{\theta(\alpha\bar{\beta})}-\alpha\bar{\beta}h_4(0).
\end{equation}
Now, equations \eqref{h4 recall} and \eqref{wsp3} together imply
\begin{equation}\label{h"(0)=}
h_5(0)= \frac{\beta}{\theta(\alpha\bar{\beta})}+\alpha\bar{\beta}\Big(\frac{\bar{\alpha}\beta\overline{\theta(\alpha\bar{\beta})} (1+\|g\|^2)}{\bar{\beta}}\Big) = \frac{\beta[1+|\alpha|^2(1+\|g\|^2)|\theta(\alpha\bar{\beta})|^2]}{\theta(\alpha\bar{\beta})}.
\end{equation}
Again, taking the inner product with $z\theta$ on both sides of \eqref{beta+(z-ab)g_1=theta h'' same} and using \eqref{g1 perp on z theta Hardy space}, we have
\begin{equation}\label{g1 theta=h'' z}
\la g_1,\theta\ra =\la h_5,z\ra.
\end{equation}
Hence the equations \eqref{g1 theta=h'' z}, \eqref{h'(0)}, and \eqref{h4 recall} together imply
\begin{equation}\label{h''(z)}
\la h_5, z\ra=\la g_1, \theta\ra=h_4(0)=-\frac{\bar{\alpha}\beta\overline{\theta(\alpha\bar{\beta})}}{\bar{\beta}}(1+\|g\|^2).
\end{equation}
Next, taking the inner product with $z^n\theta$ ($n\geq 2$) on both sides of \eqref{beta+(z-ab)g_1=theta h'' same}, it follows by \eqref{g1 perp on z theta Hardy space}, $\la h_5, z^n\ra=0$, for all $n\geq 2$.
Hence $h_5$ is a linear polynomial: $h_5(z)= h_5(0)+\la h_5,z\ra z$. Then, by equations \eqref{h"(0)=} and \eqref{h''(z)},
\[
h_5(z)=\frac{\beta}{\theta(\alpha\bar{\beta})}\big[1+|\alpha|^2(1+\|g\|^2)|\theta(\alpha\bar{\beta})|^2\big]-\frac{\bar{\alpha}\beta\overline{\theta(\alpha\bar{\beta})} (1+\|g\|^2)}{\bar{\beta}}z,
\]
and finally by \eqref{r},
\begin{equation}\label{h''(z)=a+z}
h_5(z)=\frac{\bar{\alpha}\beta\overline{\theta(\alpha\bar{\beta})}(1+\|g\|^2)}{\bar{\beta}} (r-z).
\end{equation}
This completes the proof.
\end{proof}
We are now ready to prove the main theorem of this section:
\begin{Theorem}\label{r geq 1}
Let $\clm\subseteq H^2_{\alpha,\beta}$ be a nonzero, proper, closed $S_{\alpha,\beta}$-invariant subspace. Then $\clm=[\clm\ominus S_{\alpha,\beta}\clm]_{S_{\alpha,\beta}}$ if and only if $|r|\geq 1$, where $r$ is given by \eqref{r}.
\end{Theorem}
\begin{proof}
Note by Theorem \ref{codim 1} and Theorem \ref{S1 inv 1}, $\clm\ominus S_{\alpha,\beta}\clm=\C f_2$, where $f_2$ is given by \eqref{isp4} and $\clm$ satifies the conditions \eqref{isp1}--\eqref{isp4}. Since
\[
[f_2]_{S_{\alpha,\beta}}=\overline{span}\{f_2, S_{\alpha,\beta}f_2, S_{\alpha,\beta}^2f_2, S_{\alpha,\beta}^3f_2,\ldots\}
= \overline{span}\{f_2, S_{\alpha,\beta}f_2, S_{\alpha,\beta}^2f_2-\alpha\bar{\beta}S_{\alpha,\beta}f_2,\ldots\},
\]
we can write
\begin{equation}\label{[f2]}
[f_2]_{S_{\alpha,\beta}}=\C f_2+\C S_{\alpha,\beta}f_2+[S_{\alpha,\beta}^2f_2-\alpha\bar{\beta}S_{\alpha,\beta}f_2]_{S_{\alpha,\beta}}.
\end{equation}
We now proceed for a representation of $[S_{\alpha,\beta}^2f_2-\alpha\bar{\beta}S_{\alpha,\beta}f_2]_{S_{\alpha,\beta}}$.
Note that, the equations \eqref{s12 f2 same}, \eqref{beta+(z-ab)g_1=theta h'' same}, and \eqref{h''(z)=a+z} together imply,
\begin{equation}\label{s1 squre f2 -s1 f2}
S_{\alpha,\beta}^2f_2-\alpha\bar{\beta}S_{\alpha,\beta}f_2 = z^3\theta h_5=\frac{\bar{\alpha}\beta\overline{\theta(\alpha\bar{\beta})} (1+\|g\|^2)}{\bar{\beta}}z^3\theta \big(r-z\big), \text{ and hence}
\end{equation}
\begin{equation}\label{S1f2 for [f2]}
[S_{\alpha,\beta}^2f_2-\alpha\bar{\beta}S_{\alpha,\beta}f_2]_{S_{\alpha,\beta}}=[z^3\theta(r-z)]_{S_{\alpha,\beta}}=[z^3\theta(r-z)]_S=z^3\theta[r-z]_S.
\end{equation}
Again by \eqref{[f2]} and \eqref{S1f2 for [f2]},
\begin{equation}\label{[f2] mod}
[f_2]_{S_{\alpha,\beta}}=\C f_2+\C S_{\alpha,\beta}f_2+\big[S_{\alpha,\beta}^2f_2-\alpha\bar{\beta}S_{\alpha,\beta}f_2\big]_{S_{\alpha,\beta}}=\C f_2\oplus\Big(\C S_{\alpha,\beta}f_2+z^3\theta[r-z]_S\Big),
\end{equation}
where the last equality holds as $f_2\perp S_{\alpha,\beta}\clm$ (by \eqref{isp4}). Suppose $\clm=[\clm\ominus S_{\alpha,\beta}\clm]_{S_{\alpha,\beta}}=[f_2]_{S_{\alpha,\beta}}$. Our goal is to show  $|r|\geq 1$.

Note that, $f_1\perp f_2$ (by \eqref{isp1}) and $f_2\perp S_{\alpha,\beta}\clm$ (by \eqref{isp4}) together imply $f_2\perp (f_1-S_{\alpha,\beta}f_2)$ and hence by \eqref{[f2] mod}, there exist $t_0\in\C$ and $h_1\in [r-z]_S$ such that
\begin{equation}\label{sf2-f1}
f_1-S_{\alpha,\beta}f_2= t_0 S_{\alpha,\beta}f_2+ z^3\theta h_1
\end{equation}
Recall by \eqref{s1f2=f1-theta}, $f_1-S_{\alpha,\beta}f_2=\frac{\bar{\alpha}\beta}{\bar{\beta}}\overline{\theta(\alpha\bar{\beta})}(1+\|g\|^2)z^3\theta$. Then by \eqref{isp4}, the equation \eqref{sf2-f1} further reduces to
\begin{equation}\label{t equation}
\begin{split}
\frac{\bar{\alpha}\beta}{\bar{\beta}}\overline{\theta(\alpha\bar{\beta})}(1+\|g\|^2)z^3\theta &= t_0 S_{\alpha,\beta}\Big(f_0+z^2g-\frac{\bar{\alpha}\beta}{\bar{\beta}}\overline{\theta(\alpha\bar{\beta})}(1+\|g\|^2)z^2\theta\Big)+z^3\theta h_1\\
&= t_0\Big(\alpha\bar{\beta}f_0+\beta z^2+ z^3g-\frac{\bar{\alpha}\beta}{\bar{\beta}}\overline{\theta(\alpha\bar{\beta})}(1+\|g\|^2)z^3\theta\Big)+z^3\theta h_1, \ \ (\text{ by } \eqref{eqn: S_ab matrix})
\end{split}
\end{equation}
Since L.H.S. of \eqref{t equation} is free from $z$ and $z^2$, we must have $t_0=0$ and \eqref{t equation} simplifies to
\begin{equation}\label{|h|=1}
h_1=\frac{\bar{\alpha}\beta}{\bar{\beta}}\overline{\theta(\alpha\bar{\beta})}(1+\|g\|^2).
\end{equation}
Hence, $h_1\in\big[r-z\big]_{S}$ is a nonzero constant. This implies $|r|\geq 1$.

Conversely, let $|r|\geq 1$. Then it follows by \eqref{[f2] mod},
$$[f_2]_{S_{\alpha,\beta}}=\ C f_2 \oplus \Big(\C S_{\alpha,\beta}f_2+z^3\theta H^2(\D)\Big).$$
\textsf{Claim:} $\Big(\C S_{\alpha,\beta}f_2+z^3\theta H^2(\D)\Big)=\C f_1\oplus z^3\theta H^2(\D)$.
This follows immediately by \eqref{s1f2=f1-theta}:
$$S_{\alpha,\beta}f_2=f_1-\frac{\bar{\alpha}\beta}{\bar{\beta}}\overline{\theta(\alpha\bar{\beta})}(1+\|g\|^2)z^3\theta.$$
Hence, $[f_2]_{S_{\alpha,\beta}}=\C f_2\oplus\C f_1\oplus z^3\theta H^2(\D)=\clm$.
\end{proof}
We now discuss the w.s.p. of $S_{\alpha,\beta}$-invariant subspaces $\clm$ where the associated inner functions $\theta$ are monomials. The following theorem shows that, such an $\clm$ is always generated by the corresponding wandering subspace.
\begin{Theorem}\label{z power n}
In the setting of Theorem \ref{S1 inv 1}, for $\theta(z)=z^n$, $n\geq 1$, $\clm=[\clm\ominus~ S_{\alpha,\beta}\clm]_{S_{\alpha,\beta}}$.
\end{Theorem}
\begin{proof}
Let $\theta(z)=z^n$ for some $n\in\N$. Then by Theorem \ref{S1 inv 1}, the associated function $g$ satisfies
$$(z-\alpha\bar{\beta})g=\frac{\beta}{\theta(\alpha\bar{\beta})}[z^n-(\alpha\bar{\beta})^n].$$
Hence,
$$g(z)=\frac{\beta}{(\alpha\bar{\beta})^n}[z^{n-1}+(\alpha\bar{\beta})z^{n-2}+(\alpha\bar{\beta})^2 z^{n-3}+\cdots+(\alpha\bar{\beta})^{n-2}z+(\alpha\bar{\beta})^{n-1}],$$
which further implies $\|g\|^2=\frac{|\beta|^2}{|\alpha\beta|^{2n}}\Big(\frac{1-|\alpha\beta|^{2n}}{1-|\alpha\beta|^2}\Big)$ and
\begin{equation}\label{1+norm g}
\begin{split}
1+\|g\|^2 &=1+\frac{|\beta|^2(1-|\alpha\beta|^{2n})}{|\alpha\beta|^{2n}(1-|\alpha\beta|^2)}
=\frac{|\alpha\beta|^{2n}(1-|\alpha\beta|^2)+|\beta|^2(1-|\alpha\beta|^{2n})}{|\alpha\beta|^{2n}(1-|\alpha\beta|^2)}\\
\end{split}
\end{equation}
In view of Theorem \ref{r geq 1}, we show that $|r|\geq 1$, where $r$ is given by \eqref{r}. Note that, by \eqref{r}
\begin{equation}\label{r<1}
|r|\geq 1\iff |\alpha\beta|-|\alpha|^2|\theta(\alpha\bar{\beta})|^2(1+\|g\|^2)(1-|\alpha\beta|)\geq 0
\end{equation}
Upon a simplification, the relations \eqref{r<1} and \eqref{1+norm g} together imply,
\begin{equation}\label{max alpha beta}
|r|\geq1 \iff \Big[1-\Big(|\alpha|^{2n}|\beta|^{2n-2}|\alpha\beta|(1-|\alpha\beta|)+\frac{|\alpha\beta|(1-|\alpha\beta|^{2n})}{(1+|\alpha\beta|)}\Big)\Big]\geq 0
\end{equation}
For $n=1$, $\theta(z)=z$ and (by \eqref{max alpha beta}) $|r|\geq 1$ is equivalent to
$$1-|\alpha\beta|(1-|\alpha\beta|)(1+|\alpha|^2)\geq0,$$
which is evidently true as, $\max \Big(|\alpha\beta|(1-|\alpha\beta|\Big)=\frac{1}{4}$ and  $$|\alpha\beta|(1-|\alpha\beta|)(1+|\alpha|^2)\leq\frac{1}{4}(1+|\alpha|^2)<\frac{1}{2},$$ for $0<|\alpha|<1$.
Hence for $\theta(z)=z$, $\clm=[\clm\ominus S_{\alpha,\beta}\clm]_{S_{\alpha,\beta}}$.

We now discuss the case $\theta(z)=z^n$ for $n\geq 2$. Let us denote:
\begin{equation}\label{R intro}
R:=\Big(|\alpha|^{2n}|\beta|^{2n-2}|\alpha\beta|(1-|\alpha\beta|)+\frac{|\alpha\beta|(1-|\alpha\beta|^{2n})}{(1+|\alpha\beta|)}\Big).
\end{equation}
Clearly, $R>0$ and we have
\begin{equation}\label{4a<b}
R<\frac{1}{4}|\alpha\beta|+\frac{|\alpha\beta|}{(1+|\alpha\beta|)}.
\end{equation}
\textsf{Claim:} We show
\begin{equation}\label{wsp claim}
\Big(\frac{1}{4}|\alpha\beta|+\frac{|\alpha\beta|}{(1+|\alpha\beta|)}\Big)\leq1.
\end{equation}
 Suppose on the contrary, $\Big(\frac{1}{4}|\alpha\beta|+\frac{|\alpha\beta|}{(1+|\alpha\beta|)}\Big)>1$. This is equivalent to
\begin{equation}\label{x>1}
|\alpha\beta|^2+|\alpha\beta|-4 >0.
\end{equation}
Since $|\alpha\beta|^2<|\alpha\beta|$ implies $(|\alpha\beta|^2+|\alpha\beta|)<2|\alpha\beta|<2$, we have $(|\alpha\beta|^2+|\alpha\beta|-4)<-2$, contradicting \eqref{x>1}. This settles the claim.

Now, the relations \eqref{4a<b} and \eqref{wsp claim} together imply $(1-R)>0$ and hence, $|r|\geq 1$ (by \eqref{max alpha beta} and \eqref{R intro}). Therefore, in this case also, we have $\clm=[\clm\ominus S_{\alpha,\beta}\clm]_{S_{\alpha,\beta}}$.
\end{proof}
\begin{Remark}
For $\theta(z)= z^n$ with $n=0$ i.e., for $\theta(z)=1$, one can verify that the conditions \eqref{isp1}--\eqref{isp4} will force $\clm$ to be the full space $H^2_{\alpha,\beta}$ and $f_2 \in (\clm\ominus S_{\alpha,\beta}\clm)$ will reduce to $(f_0-\frac{\bar{\alpha}\beta}{\bar{\beta}}z^2)$, which belongs to the kernel of $S_{\alpha,\beta}^*$ ( by Lemma \ref{kernel S_ab adjoint}). In this case, w.s.p. may not hold in general (see Theorem \ref{alpha less zero}).
\end{Remark}

We now consider the w.s.p. of the $S_{\alpha,\beta}$-invariant subspaces $\clm$ with associated inner functions $\theta(z)=\frac{z-a}{1-\bar{a}z}$,\ \ $0<|a|<1$. As we noted earlier, $\theta(\alpha\bar{\beta})\neq 0$, we have $a\neq \alpha\bar{\beta}$. We will follow the same technique applied for the monomials discussed above i.e., first compute the associated function $g$ \eqref{isp3} and then find a relation between $\alpha,\beta$ (and $a$) equivalent to the inequality $|r|\geq 1$ (or $|r|<1$) for further verification.

Let us set
\begin{equation}\label{t}
t=\frac{\alpha\bar{\beta}-a}{1-\bar{a}\alpha\bar{\beta}}=\theta(\alpha\bar{\beta}).
\end{equation}
Then, a simple computation will show
\begin{equation}\label{theta= t equation}
\theta(z)-\theta(\alpha\bar{\beta})=\frac{z-a}{1-\bar{a}z}-\frac{\alpha\bar{\beta}-a}{1-\bar{a}\alpha\bar{\beta}}=(1+\bar{a}t)\Big[\frac{z-\frac{a+t}{1+\bar{a}t}}{1-\bar{a}z}\Big]
,\end{equation}
and we have by \eqref{isp3},
\begin{equation}\label{z-ab g=t}
(z-\alpha\bar{\beta})g =\frac{\beta}{t}(1+\bar{a}t)\Big[\frac{z-\frac{a+t}{1+\bar{a}t}}{1-\bar{a}z}\Big].
\end{equation}
Note that, by \eqref{t}
\begin{eqnarray}
\frac{1+\bar{a}t}{t} &=& \frac{1-|a|^2}{\alpha\bar{\beta}-a}, \label{a/t}\\
\frac{a+t}{1+\bar{a}t} &=& \alpha\bar{\beta}.\label{a simp}
\end{eqnarray}


Hence, \eqref{z-ab g=t}, \eqref{a/t}, and \eqref{a simp} together imply
\begin{equation}\label{g}
(z-\alpha\bar{\beta})g(z)= \beta\Big(\frac{1-|a|^2}{\alpha\bar{\beta}-a}\Big) \frac{z-\alpha\bar{\beta}}{1-\bar{a}z}.
\end{equation}
Since $(z-\alpha\bar{\beta})\neq 0$, it follows by \eqref{g}
\begin{equation}\label{g=c/z}
g(z)\ =\ \frac{\beta(1-|a|^2)}{(\alpha\bar{\beta}-a)} \frac{1}{1-\bar{a}z}\ =\ \frac{\beta(1-|a|^2)}{(\alpha\bar{\beta}-a)} (1+\bar{a}z+\bar{a}^2z^2+\bar{a}^3z^3+\cdots).
\end{equation}
Hence, we have by \eqref{g=c/z}
\begin{equation}\label{|g|}
\begin{split}
\|g\|^2 =\frac{|\beta|^2(1-|a|^2)}{|\alpha\bar{\beta}-a|^2},
\end{split}
\end{equation}
which further implies
\begin{equation}\label{1+g}
\begin{split}
1+\|g\|^2 = \frac{|\alpha\bar{\beta}-a|^2+|\beta|^2(1-|a|^2)}{|\alpha\bar{\beta}-a|^2}
\end{split}
\end{equation}
Now substituting $t$ and $(1+\|g\|^2)$ (from \eqref{t},\eqref{1+g}) in $r$ \eqref{r}, a simple computation reveals
\begin{equation}\label{Blaschke expre}
|r|<1\iff |\alpha\beta|-\frac{|\alpha|^2}{|1-\bar{a}\alpha\bar{\beta}|^2}[|\alpha\bar{\beta}-a|^2+|\beta|^2(1-|a|^2)](1-|\alpha\beta|)< 0
\end{equation}
Note that,
\begin{eqnarray}
|\alpha\bar{\beta}-a|^2 &=& |\alpha\beta|^2+|a|^2-2\text{ Re }(\bar{a}\alpha\bar{\beta}), \ \ \label{|ab-a|}\\
|1-\bar{a}\alpha\bar{\beta}|^2 &=& 1+|a\alpha\beta|^2-2\text{ Re }(\bar{a}\alpha\bar{\beta}), \label{|1-ab|*2}
\end{eqnarray}
and hence, \eqref{|ab-a|} and \eqref{|1-ab|*2} together imply
\begin{equation}\label{|ab-a|2}
|\alpha\bar{\beta}-a|^2=|1-\bar{a}\alpha\bar{\beta}|^2-(1-|a|^2)(1-|\alpha\beta|^2)
\end{equation}
By \eqref{|ab-a|2} and $|\alpha|^2+|\beta|^2=1$, a simple computation shows
\begin{equation}\label{r within}
|\alpha\bar{\beta}-a|^2+|\beta|^2(1-|a|^2)=|1-\bar{a}\alpha\bar{\beta}|^2\Big(1-\frac{|\alpha|^4(1-|a|^2)}{|1-\bar{a}\alpha\bar{\beta}|^2}\Big).
\end{equation}
Hence, \eqref{Blaschke expre} and \eqref{r within} together imply
\begin{equation}\label{wsp inq}
|r|<1\iff |\alpha\beta|-|\alpha|^2\Big[1-\frac{|\alpha|^4(1-|a|^2)}{|1-\bar{a}\alpha\bar{\beta}|^2}\Big](1-|\alpha\beta|)<0
\end{equation}
\textsf{Claim:\ \  $0<1-\frac{|\alpha|^4(1-|a|^2)}{|1-\bar{a}\alpha\bar{\beta}|^2}<1$.} Note by \eqref{|ab-a|2} and \eqref{r within},
\begin{equation}\label{A}
1-\frac{|\alpha|^4(1-|a|^2)}{|1-\bar{a}\alpha\bar{\beta}|^2} = \frac{|\alpha\bar{\beta}-a|^2+(1-|a|^2)|\beta|^2}{|\alpha\bar{\beta}-a|^2+(1-|a|^2)(1-|\alpha\beta|^2)}:=A  \text{ (say)}
\end{equation}
Clearly $A> 0$. Again, $|\beta|^2+|\alpha\beta|^2<|\beta|^2+|\alpha|^2=1$. Then, we are done by-
\begin{align*}
\begin{aligned}
|\beta|^2 < 1-|\alpha\beta|^2 \iff & (1-|a|^2)|\beta|^2 <(1-|a|^2)(1-|\alpha\beta|^2)\ \ (\text{ as } |a|<1) \\
\iff & |\alpha\bar{\beta}-a|^2+(1-|a|^2)|\beta|^2 <|\alpha\bar{\beta}-a|^2+(1-|a|^2)(1-|\alpha\beta|^2)\\
\iff & A<1.
\end{aligned}
\end{align*}
Now \textsf{let us suppose there exists $|\alpha|,|\beta|,|a|\in (0,1)$ such that $|r|<1$ holds.} Then by \eqref{wsp inq},
\begin{equation}\label{r- right}
|\alpha\beta|-|\alpha|^2\Big[1-\frac{|\alpha|^4(1-|a|^2)}{|1-\bar{a}\alpha\bar{\beta}|^2}\Big](1-|\alpha\beta|)<0.
\end{equation}
Now the inequality \eqref{r- right}, together with the claim above implies
\begin{equation}\label{alpha-extra}
|\alpha\beta|< |\alpha|^2(1-|\alpha\beta|).
\end{equation}
Substituting $|\beta|^2=(1-|\alpha|^2)$ in \eqref{alpha-extra}, a simplification shows
\begin{equation}\label{beta < 0}
|\alpha\beta| < |\alpha|^2(1-|\alpha\beta|)\iff |\beta|^6-4|\beta|^4+5|\beta|^2-1 <0
\end{equation}
Let us set, $|\beta|^2=x$. Then, $|\beta|^6-4|\beta|^4+5|\beta|^2-1=x^3+(1-x)(4x-1)$, and hence
\begin{equation}\label{beta=x}
|\beta|^6-4|\beta|^4+5|\beta|^2-1 <0 \iff x^3+(1-x)(4x-1)<0
\end{equation}
Since $x^3$ and $(1-x)$ are positive, the inequality $x^3+(1-x)(4x-1)<0$ implies $(4x-1)<0$. Clearly, $4x-1<0$ is equivalent to $x=\frac{1}{4+y}$ for some $y>0$ and we have
\begin{equation}\label{y>0}
x^3+(1-x)(4x-1)<0\iff  y(3+y)(4+y)-1 >0
\end{equation}
Let us consider the function $f: \mathbb{R}\longrightarrow\mathbb{R}$ defined by
$$f(y)=y(3+y)(4+y)-1.$$
Note that, $f(0)=-1$ and $f(1)=19>0$. Also $f'(y)=3y^2+14y+12>0$ for all $y\geq 0$ i.e., the function is strictly increasing in $[0, \infty)$ and has exactly one root in $[0, \infty)$. It is easy to verify that the other two roots $f(y)$ are imaginary. let us denote this real root by $\gamma$ i.e., $f(\gamma)=0$, where $\gamma\in (0, 1)$. Now for $y\geq 0$,
\begin{equation}\label{y and gamma}
f(y)>0\iff y>\gamma\iff\frac{1}{4+y}<\frac{1}{4+\gamma}.
\end{equation}
Hence, the relations \eqref{y>0} and \eqref{y and gamma} together imply
\begin{equation}\label{x<gamma}
x^3+(1-x)(4x-1)<0\iff x(=\frac{1}{4+y})<\frac{1}{4+\gamma}\iff |\beta|^2 < \frac{1}{4+\gamma}.
\end{equation}
Therefore if there exist $\alpha,\beta$ and $0<|a|<1$ such that the relation \eqref{r- right} holds, then it necessarily follows,
$|\beta|^2<\frac{1}{4+\gamma}$, where $\gamma$ is the unique real root of the cubic polynomial $f(y)=y^3+7y^2+12y-1$.
Formally, this together with Theorem \ref{r geq 1} summarizes:
\begin{Theorem}\label{wsp criteria}
Let $\clm\subseteq H^2_{\alpha,\beta}$ be an $S_{\alpha,\beta}$-invariant subspace as in Theorem \ref{S1 inv 1} with the associated inner function $\theta(z)=\frac{z-a}{1-\bar{a}z}$ with $0<|a| <1$ and $a\neq \alpha\bar{\beta}$. If $\clm\neq[\clm\ominus S_{\alpha,\beta}\clm]_{S_{\alpha,\beta}}$ then $|\beta|^2<\frac{1}{4+\gamma}$, where $\gamma$ is the unique real root of the cubic polynomial $z^3+7z^2+12z-1$.
\end{Theorem}
\NI However, if $|\beta|^2<\frac{1}{4+\gamma}$, $\gamma$ being the real root of the cubic polynomial $f(y)$ above, the inequality \eqref{r- right} need not hold for all $0<|a|<1$ with $a\neq \alpha\bar{\beta}$ in general. To see this, we first note that $\gamma\in(\frac{7}{100},\frac{2}{25})$. if we choose $\beta=\frac{1}{\sqrt{5}}$ and $a=\frac{1}{2}$, then $\alpha\beta=\frac{2}{5}$ and hence, $a\neq\alpha\bar{\beta}$. Also, $|\beta|^2<\frac{1}{4+\gamma}$. Then a routine verification will show L.H.S. of \eqref{r- right} is $\big(\frac{7}{25}\big)$, strictly positive.
Towards the converse of Theorem \ref{wsp criteria}, we can say:
\textsf{If $|\beta|^2<\frac{1}{4+\gamma}$, there exists $0<|a|<1$ with $a\neq \alpha\bar{\beta}$ such that \eqref{r- right} holds.} In fact, the inequalities \eqref{beta < 0}--\eqref{x<gamma} imply
\begin{equation}\label{ab inq}
|\beta|^2<\frac{1}{4+\gamma}\iff |\alpha\beta|<|\alpha|^2(1-|\alpha\beta|) \text{ and hence},
\end{equation}
there exists an $\epsilon>0$ such that
\begin{equation}\label{epsilon}
|\alpha\beta|-|\alpha|^2(1-|\alpha\beta|)+\epsilon<0.
\end{equation}
Note that, the L.H.S. of \eqref{r- right} can be written as
\begin{equation}\label{wsplast}
|\alpha\beta|-|\alpha|^2\Big[1-\frac{|\alpha|^4(1-|a|^2)}{|1-\bar{a}\alpha\bar{\beta}|^2}\Big](1-|\alpha\beta|)=\Big(|\alpha\beta|-|\alpha|^2(1-|\alpha\beta|)\Big)+ B, \text{ where}
\end{equation}
$B:= \frac{|\alpha|^6(1-|\alpha\beta|)(1-|a|^2)}{|1-\bar{a}\alpha\bar{\beta}|^2}$. Clearly, $B$ is strictly positive and $B\longrightarrow 0$ as $|a|$ tends to $1$. Hence corresponding to the $\epsilon$ above, there exists $|a|<1$ with $a\neq\alpha\bar{\beta}$ such that
\begin{equation}\label{<e}
B<\epsilon.
\end{equation}
Hence it follows by \eqref{epsilon}--\eqref{<e},
$$|\alpha\beta|-|\alpha|^2\Big[1-\frac{|\alpha|^4(1-|a|^2)}{|1-\bar{a}\alpha\bar{\beta}|^2}\Big](1-|\alpha\beta|)<0, \text{ and therefore}$$ by \eqref{wsp inq}, $|r|<1$. As a consequence, we have the following Theorem:
\begin{Theorem}\label{not wsp}
Let $|\beta|^2<\frac{1}{4+\gamma}$, where $\gamma$ is the unique real root of the cubic polynomial $z^3+~7z^2+12z-1$. Then there exists $S_{\alpha,\beta}$-invariant subspace $\clm\subseteq H^2_{\alpha,\beta}$ with associate inner function $\theta(z)=\frac{z-a}{1-\bar{a}z}\ \ (0<|a|<1,  a\neq \alpha\bar{\beta})$, such that $\clm\neq[\clm\ominus S_{\alpha,\beta}\clm]$.
\end{Theorem}
\vspace{0.1 in}

\textbf{Acknowledgement:}
The author is thankful to Prof. E. K. Narayanan for many helpful discussions, suggestions, and corrections throughout the work. This work is supported in part by IoE-IISc Postdoctoral fellowship and NBHM-Postdoctoral fellowship of the author at the Indian Institute of Science, Bangalore, India.

\bibliographystyle{amsplain}

\end{document}